\documentclass[12pt,twoside]{article}
\usepackage{epsfig}
\usepackage{subfigure}

\usepackage{amsmath}
\usepackage{stmaryrd}
\usepackage{euscript}

\title{Particle approximation of initial data for systems of conservation laws
in two dimensions}

\author{F. Berthelin}\date{}

\begin{document}

\font\bba=msbm10
\font\bbb=msbm8
\font\bbc=msbm6
\newfam\bbfam
\textfont\bbfam=\bba
\scriptfont\bbfam=\bbb
\scriptscriptfont\bbfam=\bbc
\def\bb{\fam\bbfam\bba}

\def\N{{\bb N}}
\def\Z{{\bb Z}}
\def\r{{\bb R}}
\def\C{{\bb C}}
\def\ML{{\cal M}_t([0,+\infty[, L^\infty_{xy}(\r^2))}
\def\1{{1\hspace{-1.2mm}{\rm I}}}
\def\emptyset{\hbox{$\displaystyle/\kern -5.97pt\circ$}}
\def\diam{\mbox{{\normalsize diam}}}
\def\supp{\mbox{{\normalsize supp}}}
\def\lip{\mbox{{\normalsize Lip}}}
\def\stackunder#1#2{\mathop{#1}\limits_{#2}}
\def\stackover#1#2{\mathop{#1}\limits^{#2}}
\def\CQFD{\unskip\kern 6pt\penalty 500%
\raise -1pt\hbox{\vrule\vbox to 8pt{\hrule width 6pt\vfill\hrule}\vrule}}
\def\ess{\mathop{\mbox{\rm ess}}}
\def\limsup{\mathop{\overline{\lim}}\limits}
\def\liminf{\mathop{\underline{\lim}}\limits}
\def\dv{\mathop{\mbox{\rm div}}\nolimits}
\def\sgn{\mathop{\mbox{\rm sgn}}}
\def\iint{\mathop{\int\mkern -12mu\int}}
\def\iiint{\mathop{\int\mkern -12mu\int\mkern -12mu\int}}
\def\iiiint{\mathop{\int\mkern -12mu\int\mkern -12mu\int\mkern -12mu\int}}
\def\tq{\mathrel{;\,}}
\def\pref#1{(\ref{#1})}
\def\pr{\noindent {\bf Proof. }}
\let\dsp=\displaystyle
\let\ds=\displaystyle
\let\eps=\varepsilon
\def\convol{\mathop{*}\limits}
\def\tr{\mathop{\mbox{\rm tr}}}
\def\meas{\mathop{\rm meas}}
\def\an{\alpha_n}
\def\une{u_n^\varepsilon}
\newcommand{\integra}[1]{\int_{[0,T]} \int_\r \int_\r #1  \,d\xi \,dx \,dt}
\newcommand{\integr}[1]{\int_\r \int_\r #1  \,d\xi \,dx}
\newcommand{\integ}[1]{\int_\r  #1 \,d\xi}
\newcommand{\I}[1]{\int_0^1 #1 \,d\vartheta}
\newcommand{\intbr}[2]{\int \!\!\! \int \!\!\! \int_{\{x \in B_R\}
\cap \{#1\}} #2 \,d\xi \,dx \,dt}
\newcommand{\intO}[1]{\int_{[0,+\infty[} \iint_{\r^2}  #1  \,dx \,dy \,dt}
\newcommand{\intOb}[1]{\int_0^{+\infty}\!\!\!\int_\r  #1  }
\newcommand{\intome}[1]{\int_{[0,T]} \int_{\omega_t}  #1  \,dx \,dt}
\newcommand{\intom}[1]{\int\!\!\!\!\int_{\Omega}  #1  \,dx \,dt}
\newcommand{\intoms}[1]{\int\!\!\!\!\int_{\Omega}  #1 }
\newcommand{\intab}[1]{\int_0^\varepsilon \int_{a(t)}^{b(t)}  #1  \,dx \,dt}
\def\liu{L^{\infty}_t(0,T;L^1_x(\r))}
\def\lii{L^{\infty}_t(0,T;L^\infty_x(\r))}
\def\luu{L^1_t(0,T;L^1_x(\r))}
\def\lu{L^1(\r)}
\def\l1c{L^1_c(\r)}
\def\li{L^\infty(\r)}
\def\doi{{\cal D}([0,+\infty[\times\r)}
\def\do{{\cal D}([0,T]\times\r)}
\def\dw{{\cal D}(\Omega)}
\def\ctw{C_t([0,T];L^\infty_{w*}(\r))}
\def\tw{\mathop{{\rightharpoonup}}\limits}
\def\tow{\tw_{n \to + \infty}}
\def\m1{{\cal M}_t([0,T];L^1_x(\r))}
\def\essinf{\mathop{\ess\inf}\limits}
\def\esssup{\mathop{\ess\sup}\limits}

\def\limd{\mathop{{\lim}}\limits}
\def\tod{\mathop{{\to}}\limits}
\def\tos{\tod_{n \to + \infty}}
\def\infd{\mathop{{\inf}}\limits}
\def\supd{\mathop{{\sup}}\limits}
\def\simd{\mathop{{\sim}}\limits}
\def\ntod{\mathop{{\not\to}}\limits}
\def\cupd{\mathop{{\cup}}\limits}

\newcommand{\ent}[2]{\llbracket #1, #2 \rrbracket}

\newtheorem{Th}{Theorem}[section]
\newtheorem{Prop}[Th]{Proposition}
\newtheorem{Constr}[Th]{Construction}
\newtheorem{Lemma}[Th]{Lemma}
\newtheorem{Defin}[Th]{Definition}
\newtheorem{Cor}[Th]{Corollary}
\newtheorem{Rk}{\sl Remark}[section]
\newtheorem{Ex}{\sl Example}[section]

\def\theequation{\thesection.\arabic{equation}}
\def\thesection{\arabic{section}}
\def\thesubsection{\arabic{section}.\arabic{subsection}}
\def\thesubsubsection{\arabic{section}.\arabic{subsection}.\arabic{subsubsection}}
\newcommand\Section{%
\def\thesubsection{\arabic{section}}
\setcounter{Th}{0}
\setcounter{Rk}{0}
\setcounter{Ex}{0}
\setcounter{equation}{0}\section}
\newcommand\Subsection{%
\def\thesubsection{\arabic{section}.\arabic{subsection}}
\subsection}
\newcommand\Subsubsection{\subsubsection}

\maketitle

\bigskip

\begin{center}
        Universit\'e C\^ote d'Azur\\
        INRIA Sophia Antipolis, Project Team COFFEE\\
        Laboratoire J. A. Dieudonn\'e, UMR 7351 CNRS,\\
        Universit\'e C\^ote d'Azur, Parc Valrose,\\
        06108 Nice cedex 2, France \\
        e-mail: Florent.BERTHELIN@univ-cotedazur.fr
\end{center}
\vspace{2cm}

\begin{abstract}
In this paper, we prove particle approximations of initial data for systems of conservation laws
in two dimensions. 
This involves approaching the density but also all the additional quantities that could be verified by the model considered. We prove
that according to the hypothesis of regularity or support, the speed of convergence is of form $C/N$ or $C/N^2$.
\end{abstract}
\vspace{15mm}


\medskip

\vspace{1cm}

\newpage
\baselineskip=12pt
\parskip=0pt plus 1pt

\tableofcontents

\section{Introduction}\label{Section 1}

\subsection{Context}

In the area of systems of conservation laws, the study of approximations by particle models is important and the mathematical analysis is most of the time not advanced enough.
Some results have appeared recently for scalar conservation laws 
(\cite{DFR15}, \cite{DFFR16}, \cite{DFFR17}, \cite{RS})
and some conservation law systems (\cite{Berthelin_Goatin}).
All being placed in one dimension of space. The purpose of this paper is to provide tools for higher dimensional case studies.
In paper \cite{Berthelin_Goatin}, the first step consists in approximating the initial data by a system of particles and demonstrating in the sense of distributions the limit of this approximation towards the initial data which are the density, the momentum and other quantities considered by the system.
An initial data $\rho^0 \in L^1(\r)$ being fixed, the particles are initially placed in the following locations:
$$
\overline{x}_1^N= \sup \left\{ x\in \r \, ; \, \int_{-\infty}^{x} \rho^0(x)\,dx <l_N \right\},
$$
$$
\overline{x}_i^N= \sup 
\left\{ x\in \r \, ; \, \int_{\overline{x}_{i-1}^N}^{x} \rho^0(x)\,dx <l_N \right\},
~ \textrm{ for } i=2,\dots, N-1.
$$
where
$\ds l_N= \frac{1}{N} \int_\r \rho^0(x) \,dx$.
Having no order relation and therefore no upper bound in dimension greater than one,
this approach cannot be used for dimension two and greater so we have to proceed differently.
In \cite{Berthelin_Goatin}, the system studied is the following:
$$ 
\left\{ \begin{array}{l}
    \dsp \partial_t\rho+\partial_x(\rho v)=0,\\
    \dsp \partial_t(\rho (v+p))+\partial_x(\rho v (v+p))=0,
    \end{array} \right.
$$
and the important first step of the proof is to prove that approximated solutions related to the particles
satisfied for the initial data that
$$
\int_\r \hat\rho^N(x) \varphi(x) \,dx \tod_{N\to +\infty}  \int_\r \rho^0(x) \varphi(x) \,dx,
$$
$$
\int_\r \hat\rho^N(x) \hat v^N(x) \varphi(x) \,dx \tod_{N\to +\infty}  \int_\r \rho^0(x) v^0(x)\varphi(x) \,dx
$$
and
$$
\int_\r \hat\rho^N(x) \hat p^N(x) \varphi(x) \,dx \tod_{N\to +\infty}  \int_\r \rho^0(x) p^0(x)\varphi(x) \,dx.
$$
for all $\varphi\in\mathcal{C}^\infty_c (\r)$.
The proof makes significant use of being in one dimension of space. 
In this paper, we will present particle approximation results in two space dimensions. It can then be
extended to higher dimensions, 
the difficulty being to overcome the additional properties of $\r$ with respect to $\r^d$.
We will thus define particle approximations and prove the convergence for the density
$$
   \iint_{\r^2} \hat{\rho}^N(x,y) \varphi(x,y) \,dxdy  - \iint_{\r^2} \rho^0(x,y) \varphi(x,y) \,dxdy  \tod_{N\to+\infty} 0
$$
and for other additional quantities which are of the form
$$
   \iint_{\r^2} \hat{\rho}^N(x,y) \hat{\omega}^N(x,y) \varphi(x,y) \,dxdy  - \iint_{\r^2} \rho^0(x,y)
\omega^0(x,y) \varphi(x,y) \,dxdy \tod_{N\to+\infty} 0.
$$
For example for a multi-dimensional system like the one found in article \cite{Berthelin_Traffic_2D}, that is to say
\begin{equation} \left\{ \begin{array}{l} \label{systeme}
    \dsp \partial_t\rho+\partial_x(\rho u)+\partial_y(\rho v)=0,\\
    \dsp \partial_t(\rho (u+p))+\partial_x(\rho u(u+p))+\partial_y(\rho (u+p)v)=0, \\
    \dsp \partial_t(\rho (v+q))+\partial_x(\rho u(v+q))+\partial_y(\rho v(v+q))=0,
	\end{array} \right.
\end{equation}
the quantity $\omega$ can represent four different physical quantities: $u, v, p$ and $q$.
It is therefore interesting to be able to bundle the study of those four simultaneously.

We will demonstrate results according to the hypothesis on density,
$\rho^0 \in W^{1,\infty}(\r^2)$ or $\rho^0 \in L^1(\r^2)\cap L^\infty(\r^2)$, and with compact support or not.
Furthermore, instead of showing convergence to 0, we will present accurate estimates of how this convergence takes place. We will see in particular that according to the hypothesis of regularity or support, the speed of the convergence is of form $C/N$ or $C/N^2$.
When the initial data does not have compact support, we will obtain for any epsilon an increase by taking this parameter into account. Typically, this reflects the smallness of the density outside of a compact and is sufficient to allow us to continue this study.
We present some independent results of partial differential equations that could be studied via their approximations in order to be as general as possible.
The present work presents some similarities of approach with different works which do not focus on particle approximations but seek to approach through other objects: sticky particles \cite{BreGre}, \cite{ERS},
sticky blocks \cite{Berthelin}, \cite{BDDR}, \cite{Berthelin_Broizat}.
For these different objects, the first two-dimensional extensions 
came up in the recent past years
\cite{Berthelin_2D}, \cite{Berthelin_Traffic_2D}.

In this next two subsections, we present the particle approximation and the four general results we get.
In section 2, we establish the part of the results related to the density and in section 3, 
we prove the remaining results related to the additionnal quantities.

\subsection{Particles approximation}

Let $\rho^0 \in L^1(\r^2)$ such that $\rho^0 \geq 0$.
Let $L=(L_1, L_1',L_2,L_2')\in\r^4$ be such that
\begin{equation} \label{hypotheseL}
L_1'<L_1 \quad \textrm{ and } \quad L_2'<L_2
\end{equation}
and set $\Delta_1 = L_1-L_1'$, $\Delta_2 = L_2-L_2'$ and $C_{L}=[L_1',L_1]\times[L_2',L_2]$.
Let $N\in \N^*$.
We set
\begin{equation} \label{defpreli1}
\overline{x}_k^N = L_1' + k \frac{\Delta_1}{N}, \quad 
\overline{y}_k^N = L_2' + k \frac{\Delta_2}{N} \quad \textrm{ for any } k \in \ent{0}{N}.
\end{equation}
Notice that $\overline{x}_0^N=L_1'$, $\overline{x}_N^N=L_1$, 
$\overline{y}_0^N=L_2'$ and $\overline{y}_N^N=L_2$.
For $i \in \ent{0}{N-1}$ and $j \in \ent{0}{N-1}$, we set
\begin{equation}  \label{defpreli2}
a_{i,j}^N= \frac{N^2}{\Delta_1 \Delta_2} \int_{\overline{x}_i^N}^{\overline{x}_{i+1}^N}  \int_{\overline{y}_j^N}^{\overline{y}_{j+1}^N} \rho^0(x,y) \,dydx.
\end{equation}
We define the piecewise constant  density $\hat{\rho}^N$ by
\begin{equation} \label{densitediscrete}
\hat{\rho}^N(x,y)=\sum_{i=0}^{N-1} \sum_{j=0}^{N-1} a_{i,j}^N \1_{[\overline{x}_i^N,\overline{x}_{i+1}^N[}(x) \1_{[\overline{y}_j^N,\overline{y}_{j+1}^N[}(y).
\end{equation}
Particles approaching the density $\rho^0$ are located at points
$(\overline{x}_i^N,\overline{y}_j^N)$ and the density $\hat{\rho}^N$
is a way to make link between the particle dynamics and the macroscopic quantity.

Let $\rho^0 \in L^1(\r^2)$ such that $\rho^0 \geq 0$.
Let $L=(L_1, L_1',L_2,L_2')\in\r^4$ be such that (\ref{hypotheseL}) 
and set $\Delta_1 = L_1-L_1'$, $\Delta_2 = L_2-L_2'$ and $C_{L}=[L_1',L_1]\times[L_2',L_2]$.
Let $N\in \N^*$.
We set (\ref{defpreli1}) and 
for $i \in \ent{0}{N-1}$ and $j \in \ent{0}{N-1}$, we set
(\ref{defpreli2}).
We define the piecewise constant  density $\hat{\rho}^N$ by
(\ref{densitediscrete}).

Let $\omega^0 \in L^\infty(\r^2)$. Then $\omega^0 \in L^1_{loc}(\r^2)$ and
we set
\begin{equation}\label{initialdiscretquantity}
W_{i,j}^N= \frac{N^2}{\Delta_1 \Delta_2} \int_{\overline{x}_i^N}^{\overline{x}_{i+1}^N}  \int_{\overline{y}_j^N}^{\overline{y}_{j+1}^N} \omega^0(x,y) \,dydx,
 \quad \textrm{ for } i, j \in \ent{0}{N-1}
\end{equation}
and we define the piecewise constant quantity $\hat{\rho}^N \hat{\omega}^N$ by
\begin{equation} \label{quantitydiscrete}
(\hat{\rho}^N \hat{\omega}^N)(x,y)=\sum_{i=0}^{N-1} \sum_{j=0}^{N-1} a_{i,j}^N W_{i,j}^N
 \1_{[\overline{x}_i^N,\overline{x}_{i+1}^N[}(x) \1_{[\overline{y}_j^N,\overline{y}_{j+1}^N[}(y).
\end{equation}


From that starting point, 
in the paper, 
we will consider two cases of locations for the particles according to the 
assumptions on the initial data.

\begin{enumerate} 
\item {\bf Case 1 of initial data.}
In the case where $\rho^0 \in L^1(\r^2)$ has a compact support included
in some $C_L$, we choose this $L$ for the location area of
particles.
\item {\bf Case 2 of initial data.}
In the case where $\rho^0 \in L^1(\r^2)$ doesn't not a compact support,
We will fix the location area of the particles by removing a mass 
smaller than an epsilon from this area.
Indeed, since $\rho^0 \in L^1(\r^2)$ and such that $\rho^0 \geq 0$,
the function $L\mapsto \ds
\iint_{\r^2\setminus [-L,L]^2} \rho^0(x,y) \,dxdy$
is decreasing and tends to 0, then for any $\eps>0$, there exists $L>0$ such that
\begin{equation} \label{hypotheseepsL}
\iint_{\r^2\setminus [-L,L]^2} \rho^0(x,y) \,dxdy \leq \eps.
\end{equation}
Then we set $L_1=L_2=L$, $L_1'=L_2'=-L$ which gives $C_{L}=[-L,L]^2$.
Notice that in this case $\hat{\rho}^{N}$ depends on $\eps$ via $L$
thus we note $\hat{\rho}^{N,\eps}$ instead  to keep it in mind.
\end{enumerate}

\subsection{Results}

We prove the following approximation results according to the regularity of $\rho^0$, namely 
$W^{1,\infty}(\r^2)$ or $\rho^0 \in L^{1}(\r^2)\cap L^{\infty}(\r^2)$, 
and according to its support, compact or not.

\begin{Th}[Result for $\rho^0 \in W^{1,\infty}(\r^2)$ with a compact support] \label{TH1}
Let $\rho^0 \in L^1(\r^2)$ such that $\rho^0 \geq 0$ 
with a compact support included in a $C_{L}$
where $L=(L_1, L_1',L_2,L_2')\in\r^4$ be such that
(\ref{hypotheseL}).
We assume furthermore that $\partial_x \rho^0, \partial_y \rho^0 \in L^\infty(\r^2)$.
For any $N\in \N^*$, we consider
$\hat{\rho}^N$ defined by (\ref{densitediscrete}) with 
(\ref{defpreli1})-(\ref{defpreli2}).
Then for any function $\varphi \in C^\infty_c(\r^2)$,
we have
\begin{equation} \label{Estimation1a}
  \left| \iint_{\r^2}  \hat{\rho}^{N}(x,y) \varphi(x,y) \,dxdy  - \iint_{\r^2} \rho^0(x,y) \varphi(x,y) \,dxdy \right| 
\leq    \frac{C_{1,2}(\varphi)}{N^2},
\end{equation}
with
\begin{equation} \label{C12}
C_{1,2}(\varphi) =   \max(\Delta_1,\Delta_2)^4 (\|\partial_x \rho^0\|_\infty +
\|\partial_y \rho^0\|_\infty)
 ( \|\partial_x \varphi\|_\infty+ \|\partial_y \varphi\|_\infty).
\end{equation}
Let $\omega^0 \in L^\infty(\r^2)$ such that $\partial_x \omega^0, \partial_y \omega^0 \in L^\infty(\r^2)$.
Let $\hat{\rho}^N \hat{\omega}^N$ defined by (\ref{quantitydiscrete}) with (\ref{initialdiscretquantity}).
Then for any function $\varphi \in C^\infty_c(\r^2)$,
we have
\begin{equation} \label{Estimation1b}
  \left|  \iint_{\r^2} \hat{\rho}^N(x,y) \hat{\omega}^N(x,y) \varphi(x,y) \,dxdy  - \iint_{\r^2} \rho^0(x,y)
\omega^0(x,y) \varphi(x,y) \,dxdy  \right| 
\leq     \frac{K_{1,2}(\varphi)}{N} +\frac{C_{1,2}(\varphi) \|\omega^0\|_\infty}{N^2},
\end{equation}
with 
\begin{equation} \label{K12b}
K_{1,2}(\varphi) =
\max(\Delta_1,\Delta_2) (\|\partial_x \omega^0\|_\infty +
\|\partial_y \omega^0\|_\infty) \|\varphi\|_\infty \| \rho^0 \|_1.
\end{equation}
\end{Th}

\begin{Th}[Result for $\rho^0 \in W^{1,\infty}(\r^2)$]  \label{TH2}
Let $\rho^0 \in L^1(\r^2)$ and such that $\rho^0 \geq 0$.
Let $\eps >0$.
Consider $L>0$ such that (\ref{hypotheseepsL}).
We assume furthermore that $\partial_x \rho^0, \partial_y \rho^0 \in L^\infty(\r^2)$.
For any $N\in \N^*$, we consider
$\hat{\rho}^{N,\eps}$ defined by (\ref{densitediscrete}) with 
(\ref{defpreli1})-(\ref{defpreli2}).
Then for any function $\varphi \in C^\infty_c(\r^2)$,
we have
\begin{equation} \label{Estimation2a}
\left|  \iint_{\r^2} \hat{\rho}^{N,\eps}(x,y) \varphi(x,y) \,dxdy  - \iint_{\r^2} \rho^0(x,y) \varphi(x,y) \,dxdy  \right| \leq   \eps  \|\varphi \|_\infty  + \frac{C_\eps(\varphi)}{N^2},
\end{equation}
with
\begin{equation} \label{Ceps}
C_\eps(\varphi) = 16 L^4 (\|\partial_x \rho^0\|_\infty +
\|\partial_y \rho^0\|_\infty)
 ( \|\partial_x \varphi\|_\infty+ \|\partial_y \varphi\|_\infty).
\end{equation}
Let $\omega^0 \in L^\infty(\r^2)$ such that $\partial_x \omega^0, \partial_y \omega^0 \in L^\infty(\r^2)$.
Let $\hat{\rho}^{N,\eps} \hat{\omega}^{N,\eps}$ defined by (\ref{quantitydiscrete}) with (\ref{initialdiscretquantity}).
Then for any function $\varphi \in C^\infty_c(\r^2)$,
we have
\begin{align} \nonumber
& \left| \iint_{\r^2} \hat{\rho}^{N,\eps}(x,y) \hat{\omega}^{N,\eps}(x,y) \varphi(x,y) \,dxdy  -
 \iint_{\r^2} \rho^0(x,y) \omega^0(x,y) \varphi(x,y) \,dxdy \right| \\
& \qquad\qquad\qquad \leq   \eps  \|\varphi \|_\infty \|\omega^0\|_\infty + \frac{C_\eps(\varphi)\|\omega^0\|_\infty}{N^2} + \frac{K_\eps(\varphi)}{N},  \label{Estimation2b}
\end{align}
with 
\begin{equation} \label{Keps}
K_\eps(\varphi) = 2 L (\|\partial_x \omega^0\|_\infty +
\|\partial_y \omega^0\|_\infty) \|\varphi\|_\infty \|\rho^0\|_1.
\end{equation}
\end{Th}

\begin{Th}[Result for $\rho^0 \in L^{1}(\r^2)\cap L^{\infty}(\r^2)$ with a compact support] \label{TH3}
Let $\rho^0 \in L^1(\r^2)\cap L^{\infty}(\r^2)$ such that $\rho^0 \geq 0$ 
with a compact support included in a $C_{L}$
where $L=(L_1, L_1',L_2,L_2')\in\r^4$ be such that
(\ref{hypotheseL}).
For any $N\in \N^*$, we consider
$\hat{\rho}^N$ defined by (\ref{densitediscrete}) with 
(\ref{defpreli1})-(\ref{defpreli2}).
Then for any function $\varphi \in C^\infty_c(\r^2)$,
we have
\begin{equation} \label{Estimation3a}
    \left|  \iint_{\r^2} \hat{\rho}^N(x,y) \varphi(x,y) \,dxdy  - \iint_{\r^2} \rho^0(x,y) \varphi(x,y) \,dxdy  \right| \leq  \frac{D_{1,2}(\varphi)}{N}.
\end{equation}
with
\begin{equation} \label{D12}
{D}_{1,2}(\varphi) =  \max(\Delta_1,\Delta_2) 
 ( \|\partial_x \varphi\|_\infty+ \|\partial_y \varphi\|_\infty)
 (\|\rho^0\|_1 + \|\rho^0\|_\infty \Delta_1\Delta_2).
\end{equation}
Let $\omega^0 \in L^\infty(\r^2)$ such that $\partial_x \omega^0, \partial_y \omega^0 \in L^\infty(\r^2)$.
Let $\hat{\rho}^N \hat{\omega}^N$ defined by (\ref{quantitydiscrete}) with (\ref{initialdiscretquantity}).
Then for any function $\varphi \in C^\infty_c(\r^2)$,
we have
\begin{equation} \label{Estimation3b}
  \left|  \iint_{\r^2} \hat{\rho}^N(x,y) \hat{\omega}^N(x,y) \varphi(x,y) \,dxdy  - \iint_{\r^2} \rho^0(x,y)
\omega^0(x,y) \varphi(x,y) \,dxdy  \right| 
\leq     \frac{K_{1,2}(\varphi)+D_{1,2}(\varphi) \|\omega^0\|_\infty}{N},
\end{equation}
\end{Th}

\begin{Th}[Result for $\rho^0 \in L^{1}(\r^2)\cap L^{\infty}(\r^2)$]  \label{TH4}
Let $\rho^0 \in L^1(\r^2)\cap L^{\infty}(\r^2)$ and such that $\rho^0 \geq 0$.
Let $\eps >0$.
Consider $L>0$ such that (\ref{hypotheseepsL}).
For any $N\in \N^*$, we consider
$\hat{\rho}^{N,\eps}$ defined by (\ref{densitediscrete}) with 
(\ref{defpreli1})-(\ref{defpreli2}).
Then for any function $\varphi \in C^\infty_c(\r^2)$,
we have
\begin{equation} \label{Estimation4a}
\left|  \iint_{\r^2} \hat{\rho}^{N,\eps}(x,y) \varphi(x,y) \,dxdy  - \iint_{\r^2} \rho^0(x,y) \varphi(x,y) \,dxdy  \right| \leq \eps \|  \varphi \|_\infty  + \frac{D_\eps(\varphi)}{N},
\end{equation}
with
\begin{equation} \label{Deps}
D_\eps(\varphi) = 2L ( \|\partial_x \varphi\|_\infty+ \|\partial_y \varphi\|_\infty) 
( \|\rho^0\|_1  +    4L^2  \|\rho^0\|_\infty).
\end{equation}
Let $\omega^0 \in L^\infty(\r^2)$ such that $\partial_x \omega^0, \partial_y \omega^0 \in L^\infty(\r^2)$.
Let $\hat{\rho}^{N,\eps} \hat{\omega}^{N,\eps}$ defined by (\ref{quantitydiscrete}) with (\ref{initialdiscretquantity}).
Then for any function $\varphi \in C^\infty_c(\r^2)$,
we have
\begin{align} \nonumber
& \left| \iint_{\r^2} \hat{\rho}^{N,\eps}(x,y) \hat{\omega}^{N,\eps}(x,y) \varphi(x,y) \,dxdy  -
 \iint_{\r^2} \rho^0(x,y) \omega^0(x,y) \varphi(x,y) \,dxdy \right| \\
& \qquad\qquad\qquad \leq    \eps  \|\varphi \|_\infty \|\omega^0\|_\infty + \frac{K_\eps(\varphi)+D_\eps(\varphi)\|\omega^0\|_\infty}{N}.
 \label{Estimation4b}
\end{align}
\end{Th}

\begin{Rk}
We have stated the result with $\varphi \in C^\infty_c(\r^2)$ but in fact 
$\varphi \in C^\infty_c(\r^2)\cap L^\infty(\r^2)$ such that 
$\partial_x \phi, \partial_y \phi \in L^\infty(\r^2)$ is enough as we will see in the proofs.
\end{Rk}

\section{Approximation of the density of particles}

\subsection{General estimates}  \label{GR}

We start by presenting properties valid whatever are the regularity and support of density $\rho^0$.

\begin{Lemma} \label{lemmaapprox1}
Let $\rho^0 \in L^1(\r^2)$ such that $\rho^0 \geq 0$.
Let $L=(L_1, L_1',L_2,L_2')\in\r^4$ be such that
(\ref{hypotheseL})
and set $\Delta_1 = L_1-L_1'$, $\Delta_2 = L_2-L_2'$ and $C_{L}=[L_1',L_1]\times[L_2',L_2]$.
Let $N\in \N^*$ and
$\hat{\rho}^N$ defined by (\ref{densitediscrete}) with 
(\ref{defpreli1})-(\ref{defpreli2}).
Then for any function $\varphi \in C(\r^2)\cap L^\infty(\r^2)$,
we have
\begin{align*}
&   \iint_{\r^2} \hat{\rho}^N(x,y) \varphi(x,y) \,dxdy  - \iint_{\r^2} \rho^0(x,y) \varphi(x,y) \,dxdy  \\
= & - \iint_{\r^2\setminus C_{L}} \rho^0(x,y) \varphi(x,y) \,dxdy +  \sum_{i=0}^{N-1} \sum_{j=0}^{N-1} \int_{\overline{x}_i^N}^{\overline{x}_{i+1}^N}
\int_{\overline{y}_j^N}^{\overline{y}_{j+1}^N} (a_{i,j}^N-\rho^0(x,y)) ( \varphi(x,y) - \varphi(\overline{x}_i^N,\overline{y}_j^N))\,  \,dydx.
\end{align*}
\end{Lemma}
\pr 
Since $\varphi \in L^\infty(\r^2)$ and
$\rho^0 \in L^1(\r^2)$, then
$\ds \iint_{\r^2} \rho^0(x,y) \varphi(x,y) \,dxdy$ exists.
Notice that with the definition of $\hat{\rho}^N$, we have $\hat{\rho}^N=0$ on $\r^2\setminus C_{L}$, thus
we  have in one hand
\begin{align*}
 \iint_{\r^2} \hat{\rho}^N(x,y) \varphi(x,y) \,dxdy  
= & 
\iint_{C_{L}} \hat{\rho}^N(x,y) \varphi(x,y) \,dxdy \\  
= & 
\sum_{i=0}^{N-1} \sum_{j=0}^{N-1} 
\int_{\overline{x}_i^N}^{\overline{x}_{i+1}^N} \int_{\overline{y}_j^N}^{\overline{y}_{j+1}^N}  a_{i,j}^N  \varphi(x,y) \, \,dydx \\
= &  
\sum_{i=0}^{N-1} \sum_{j=0}^{N-1} 
\int_{\overline{x}_i^N}^{\overline{x}_{i+1}^N} \int_{\overline{y}_j^N}^{\overline{y}_{j+1}^N}  a_{i,j}^N  \varphi(\overline{x}_i^N,\overline{y}_j^N) \,  \,dydx\\
&  +\sum_{i=0}^{N-1} \sum_{j=0}^{N-1} a_{i,j}^N
\int_{\overline{x}_i^N}^{\overline{x}_{i+1}^N} \int_{\overline{y}_j^N}^{\overline{y}_{j+1}^N}  
 ( \varphi(x,y) - \varphi(\overline{x}_i^N,\overline{y}_j^N)) \,  \,dydx.
\end{align*}
On the other hand, we have
\begin{align*}
 \iint_{\r^2} \rho^0(x,y) \varphi(x,y) \,dxdy  
= & \iint_{\r^2\setminus C_{L}} \rho^0(x,y) \varphi(x,y) \,dxdy + \iint_{C_L} \rho^0(x,y) \varphi(x,y) \,dxdy\\
= &   \iint_{\r^2\setminus C_{L}} \rho^0(x,y) \varphi(x,y) \,dxdy + \sum_{i=0}^{N-1} \sum_{j=0}^{N-1} 
\int_{\overline{x}_i^N}^{\overline{x}_{i+1}^N} \int_{\overline{y}_j^N}^{\overline{y}_{j+1}^N} \rho^0(x,y)  \varphi(x,y) \, \,dydx,
\end{align*}
as a consequence we write
\begin{align*}
 \iint_{\r^2} \rho^0(x,y) \varphi(x,y) \,dxdy  
= &    \iint_{\r^2\setminus C_{L}} \rho^0(x,y) \varphi(x,y) \,dxdy + \sum_{i=0}^{N-1} \sum_{j=0}^{N-1} 
\int_{\overline{x}_i^N}^{\overline{x}_{i+1}^N} \int_{\overline{y}_j^N}^{\overline{y}_{j+1}^N} \rho^0(x,y)  \varphi(\overline{x}_i^N,\overline{y}_j^N)\,  \,dydx \\
&  + \sum_{i=0}^{N-1} \sum_{j=0}^{N-1} \int_{\overline{x}_i^N}^{\overline{x}_{i+1}^N}
\int_{\overline{y}_j^N}^{\overline{y}_{j+1}^N} \rho^0(x,y) (  \varphi(x,y) - \varphi(\overline{x}_i^N,\overline{y}_j^N))\,  \,dydx.
\end{align*}
We notice that
\begin{align*}
 \sum_{i=0}^{N-1} \sum_{j=0}^{N-1} 
\int_{\overline{x}_i^N}^{\overline{x}_{i+1}^N} \int_{\overline{y}_j^N}^{\overline{y}_{j+1}^N}  a_{i,j}^N   \varphi(\overline{x}_i^N,\overline{y}_j^N) \,  \,dydx 
= &   \sum_{i=0}^{N-1} \sum_{j=0}^{N-1} a_{i,j}^N \varphi(\overline{x}_i^N,\overline{y}_j^N)
\int_{\overline{x}_i^N}^{\overline{x}_{i+1}^N} \int_{\overline{y}_j^N}^{\overline{y}_{j+1}^N}  \,  \,dydx \\
= &   \sum_{i=0}^{N-1} \sum_{j=0}^{N-1} a_{i,j}^N  \varphi(\overline{x}_i^N,\overline{y}_j^N) \frac{\Delta_1 \Delta_2}{N^2} \\
= &  \sum_{i=0}^{N-1} \sum_{j=0}^{N-1} \varphi(\overline{x}_i^N,\overline{y}_j^N) 
 \int_{\overline{x}_i^N}^{\overline{x}_{i+1}^N} \int_{\overline{y}_j^N}^{\overline{y}_{j+1}^N} \rho^0(x,y)  \, \,dydx \\
= &   \sum_{i=0}^{N-1} \sum_{j=0}^{N-1} 
  \int_{\overline{x}_i^N}^{\overline{x}_{i+1}^N}  \int_{\overline{y}_j^N}^{\overline{y}_{j+1}^N}  \rho^0(x,y)  \varphi(\overline{x}_i^N,\overline{y}_j^N)
  \,  \,dydx.
\end{align*}
Then putting all this together we get
\begin{align*}
&   \iint_{\r^2} \hat{\rho}^N(x,y) \varphi(x,y) \,dxdy  - \iint_{\r^2} \rho^0(x,y) \varphi(x,y) \,dxdy  \\
= & - \iint_{\r^2\setminus C_{L}} \rho^0(x,y) \varphi(x,y) \,dxdy +  \sum_{i=0}^{N-1} \sum_{j=0}^{N-1} \int_{\overline{x}_i^N}^{\overline{x}_{i+1}^N}
\int_{\overline{y}_j^N}^{\overline{y}_{j+1}^N} (a_{i,j}^N-\rho^0(x,y)) ( \varphi(x,y) - \varphi(\overline{x}_i^N,\overline{y}_j^N))\,  \,dydx.     \CQFD
\end{align*}
Notice now the following estimate for test functions.
\begin{Lemma}  \label{lemmaapprox3preli}
Let $N\in \N^*$ and
$(\overline{x}_i^N)_{i \in \ent{0}{N-1}}$, $(\overline{y}_j^N)_{j \in \ent{0}{N-1}}$
defined by  (\ref{defpreli1}).
For any function $\varphi \in C(\r^2)\cap L^\infty(\r^2)$ such that $\partial_x \varphi, \partial_y \varphi \in L^\infty(\r^2)$,
we have
$$
  \sum_{i=0}^{N-1} \sum_{j=0}^{N-1}  \, 
\int_{\overline{x}_i^N}^{\overline{x}_{i+1}^N} \int_{\overline{y}_j^N}^{\overline{y}_{j+1}^N}  
  \left| \varphi(x,y) - \varphi(\overline{x}_i^N,\overline{y}_j^N) \right|  \,dydx
\leq      ( \|\partial_x \varphi\|_\infty+ \|\partial_y \varphi\|_\infty) 
\frac{\Delta_1 \Delta_2 \max(\Delta_1,\Delta_2)}{N}.  
$$
\end{Lemma}
\pr 
Since 
\begin{align*}
  \left| \varphi(x,y) - \varphi(\overline{x}_i^N,\overline{y}_j^N) \right| 
  \leq &   \left| \varphi(x,y) - \varphi(x,\overline{y}_j^N) \right| 
    +    \left|  \varphi(x,\overline{y}_j^N) -  \varphi(\overline{x}_i^N,\overline{y}_j^N)\right| \\
    \leq &  \|\partial_y \varphi\|_\infty (y-\overline{y}_j^N) +  \|\partial_y \varphi\|_\infty (x-\overline{x}_i^N)
\end{align*}    
for any $i, j \in \ent0{N-1}$ and any $x\in [\overline{x}_i^N ,\overline{x}_{i+1}^N]$,
$y \in [\overline{y}_j^N , \overline{y}_{j+1}^N]$,
we have
\begin{align*}
&    \sum_{i=0}^{N-1} \sum_{j=0}^{N-1}  \, 
\int_{\overline{x}_i^N}^{\overline{x}_{i+1}^N} \int_{\overline{y}_j^N}^{\overline{y}_{j+1}^N}  
  \left| \varphi(x,y) - \varphi(\overline{x}_i^N,\overline{y}_j^N) \right|  \,dydx\\
  \leq &  \sum_{i=0}^{N-1} \sum_{j=0}^{N-1} \, 
\int_{\overline{x}_i^N}^{\overline{x}_{i+1}^N} \left( \int_{\overline{y}_j^N}^{\overline{y}_{j+1}^N}  
   \|\partial_y \varphi\|_\infty (y-\overline{y}_j^N) \,dy  \right)  \,dx 
  +  \sum_{i=0}^{N-1} \sum_{j=0}^{N-1}  \, 
 \int_{\overline{y}_j^N}^{\overline{y}_{j+1}^N}  \left( \int_{\overline{x}_i^N}^{\overline{x}_{i+1}^N}
  \|\partial_y \varphi\|_\infty (x-\overline{x}_i^N)   \,dx \right)  \, dy.
\end{align*}
Thus we obtain
\begin{align*}
&  \sum_{i=0}^{N-1} \sum_{j=0}^{N-1}  \, 
\int_{\overline{x}_i^N}^{\overline{x}_{i+1}^N} \int_{\overline{y}_j^N}^{\overline{y}_{j+1}^N}  
  \left| \varphi(x,y) - \varphi(\overline{x}_i^N,\overline{y}_j^N) \right|  \,dydx \\
\leq &  \|\partial_y \varphi\|_\infty \sum_{i=0}^{N-1} \sum_{j=0}^{N-1}  \,
\int_{\overline{x}_i^N}^{\overline{x}_{i+1}^N} \left[ \frac12 (y-\overline{y}_j^N )^2 \right]_{\overline{y}_j^N}^{\overline{y}_{j+1}^N}   \,dx 
  + \|\partial_x \varphi\|_\infty \sum_{i=0}^{N-1} \sum_{j=0}^{N-1}  \,
 \int_{\overline{y}_j^N}^{\overline{y}_{j+1}^N}  \left[ \frac12 (x-\overline{x}_i^N )^2 \right]_{\overline{x}_i^N}^{\overline{x}_{i+1}^N} \, dy \\
 \leq &  \|\partial_y \varphi\|_\infty \sum_{i=0}^{N-1} \sum_{j=0}^{N-1}  \, 
\frac{\Delta_1}{N}   \frac{\Delta_2^2}{N^2}  
  +  \|\partial_x \varphi\|_\infty \sum_{i=0}^{N-1} \sum_{j=0}^{N-1}  \, \frac{\Delta_2}{N}   \frac{\Delta_1^2}{N^2}  \\
 \leq &   ( \|\partial_x \varphi\|_\infty+ \|\partial_y \varphi\|_\infty) 
 \frac{\Delta_1 \Delta_2 \max(\Delta_1,\Delta_2)}{N}.    \CQFD
\end{align*}

\subsection{Case where $\rho^0 \in W^{1,\infty}(\r^2)$}

We now move on to the estimate for the density approximation in the case of regularity 
$\rho^0 \in W^{1,\infty}(\r^2)$.
In this case, we can first estimate the difference between $a_{i,j}^N$
and $\rho^0(x,y)$ for $(x,y)\in [\overline{x}_i^N, \overline{x}_{i+1}^N]\times [\overline{y}_j^N,\overline{y}_{j+1}^N]$.
\begin{Lemma} \label{lemmaapprox1suppl}
Let $\rho^0 \in L^1(\r^2)$ such that $\rho^0 \geq 0$.
We assume furthermore that $\partial_x \rho^0, \partial_y \rho^0 \in L^\infty(\r^2)$.
Let $L=(L_1, L_1',L_2,L_2')\in\r^4$ be such that
(\ref{hypotheseL})
and set $\Delta_1 = L_1-L_1'$, $\Delta_2 = L_2-L_2'$ and $C_{L}=[L_1',L_1]\times[L_2',L_2]$.
Let $N\in \N^*$ and
$\hat{\rho}^N$ defined by (\ref{densitediscrete}) with 
(\ref{defpreli1})-(\ref{defpreli2}).
Then we have
\begin{align*}
\left|a_{i,j}^N-\rho^0(x,y)\right|  \leq  \frac{\max(\Delta_1,\Delta_2)}{N} (\|\partial_x \rho^0\|_\infty +
\|\partial_y \rho^0\|_\infty),
\end{align*}
for $(x,y)\in [\overline{x}_i^N, \overline{x}_{i+1}^N]\times [\overline{y}_j^N,\overline{y}_{j+1}^N]$.
\end{Lemma}
\pr 
We have
\begin{align*}
a_{i,j}^N-\rho^0(x,y) & = \frac{N^2}{\Delta_1 \Delta_2} \int_{\overline{x}_i^N}^{\overline{x}_{i+1}^N}  \int_{\overline{y}_j^N}^{\overline{y}_{j+1}^N} \rho^0(\tilde x,\tilde y) 
\,d\tilde yd\tilde x - \rho^0(x,y) \\
& =  \frac{N^2}{\Delta_1 \Delta_2} \int_{\overline{x}_i^N}^{\overline{x}_{i+1}^N}  \int_{\overline{y}_j^N}^{\overline{y}_{j+1}^N} (\rho^0(\tilde x,\tilde y) - \rho^0(x,y) )
\,d\tilde yd\tilde x.
\end{align*}
Then we get,
for any $x\in [\overline{x}_i^N, \overline{x}_{i+1}^N]$ and any $y \in 
[\overline{y}_j^N,\overline{y}_{j+1}^N]$,
\begin{align*}
\left|a_{i,j}^N-\rho^0(x,y)\right| & = \left|\frac{N^2}{\Delta_1 \Delta_2} \int_{\overline{x}_i^N}^{\overline{x}_{i+1}^N}  \int_{\overline{y}_j^N}^{\overline{y}_{j+1}^N} \rho^0(\tilde x,\tilde y) 
\,d\tilde yd\tilde x - \rho^0(x,y)\right| \\
& = \left| \frac{N^2}{\Delta_1 \Delta_2} \int_{\overline{x}_i^N}^{\overline{x}_{i+1}^N}  \int_{\overline{y}_j^N}^{\overline{y}_{j+1}^N} (\rho^0(\tilde x,\tilde y) - \rho^0(x,y))
\,d\tilde yd\tilde x \right| \\
& = \left| \frac{N^2}{\Delta_1 \Delta_2} \int_{\overline{x}_i^N}^{\overline{x}_{i+1}^N}  \int_{\overline{y}_j^N}^{\overline{y}_{j+1}^N} (\rho^0(\tilde x,\tilde y)- \rho^0(x,\tilde y)
+ \rho^0(x,\tilde y) - \rho^0(x,y))
\,d\tilde yd\tilde x \right| \\
& \leq  \frac{N^2}{\Delta_1 \Delta_2} \int_{\overline{x}_i^N}^{\overline{x}_{i+1}^N}  \int_{\overline{y}_j^N}^{\overline{y}_{j+1}^N} \left|\rho^0(\tilde x,\tilde y)- \rho^0(x,\tilde y)
+ \rho^0(x,\tilde y) - \rho^0(x,y) \right|
\,d\tilde yd\tilde x \\
& \leq  \frac{N^2}{\Delta_1 \Delta_2} \int_{\overline{x}_i^N}^{\overline{x}_{i+1}^N}  \int_{\overline{y}_j^N}^{\overline{y}_{j+1}^N} (\|\partial_x \rho^0\|_\infty |\tilde x -x|
+ \|\partial_y \rho^0\|_\infty |\tilde y - y| )
\,d\tilde yd\tilde x \\
& \leq  \frac{N^2}{\Delta_1 \Delta_2} \frac{\Delta_1 \Delta_2}{N^2} \|\partial_x \rho^0\|_\infty 
\int_{\overline{y}_j^N}^{\overline{y}_{j+1}^N} \,d\tilde y
+  \frac{N^2}{\Delta_1 \Delta_2} \frac{\Delta_1 \Delta_2}{N^2} \|\partial_x \rho^0\|_\infty 
\int_{\overline{x}_i^N}^{\overline{x}_{i+1}^N}  \,d\tilde x \\
& \leq  \frac{\max(\Delta_1,\Delta_2)}{N} (\|\partial_x \rho^0\|_\infty +
\|\partial_y \rho^0\|_\infty).\CQFD
\end{align*}
This estimate allows to apply the general Lemmas of section \ref{GR} and to get the following result.
\begin{Prop} \label{propapprox1}
Let $\rho^0 \in L^1(\r^2)$ such that $\rho^0 \geq 0$.
We assume furthermore that $\partial_x \rho^0, \partial_y \rho^0 \in L^\infty(\r^2)$.
Let $L=(L_1, L_1',L_2,L_2')\in\r^4$ be such that
(\ref{hypotheseL})
and set $\Delta_1 = L_1-L_1'$, $\Delta_2 = L_2-L_2'$ and $C_{L}=[L_1',L_1]\times[L_2',L_2]$.
Let $N\in \N^*$ and
$\hat{\rho}^N$ defined by (\ref{densitediscrete}) with 
(\ref{defpreli1})-(\ref{defpreli2}).
Then for any function $\varphi \in C(\r^2)\cap L^\infty(\r^2)$ such that $\partial_x \varphi, \partial_y \varphi \in L^\infty(\r^2)$,
we have
\begin{align*}
&  \left| \iint_{\r^2}  \hat{\rho}^{N}(x,y) \varphi(x,y) \,dxdy  - \iint_{\r^2} \rho^0(x,y) \varphi(x,y) \,dxdy \right| \\
\leq & \left| \iint_{\r^2\setminus C_{L}} \rho^0(x,y) \varphi(x,y) \,dxdy\right|
 +  \frac{\max(\Delta_1,\Delta_2)^4}{N^2} (\|\partial_x \rho^0\|_\infty +
\|\partial_y \rho^0\|_\infty)
 ( \|\partial_x \varphi\|_\infty+ \|\partial_y \varphi\|_\infty).
\end{align*}
\end{Prop}
\pr 
Thanks to Lemma \ref{lemmaapprox1}, we have
\begin{align*}
&   \iint_{\r^2}  \hat{\rho}^{N,\eps}(x,y) \varphi(x,y) \,dxdy  - \iint_{\r^2} \rho^0(x,y) \varphi(x,y) \,dxdy  \\
= & - \iint_{\r^2\setminus C_{L}} \rho^0(x,y) \varphi(x,y) \,dxdy +  \sum_{i=0}^{N-1} \sum_{j=0}^{N-1} \int_{\overline{x}_i^N}^{\overline{x}_{i+1}^N}
\int_{\overline{y}_j^N}^{\overline{y}_{j+1}^N} (a_{i,j}^N-\rho^0(x,y)) ( \varphi(x,y) - \varphi(\overline{x}_i^N,\overline{y}_j^N))\,  \,dydx.
\end{align*}
By using Lemma \ref{lemmaapprox1suppl}, we have
$$
\left|a_{i,j}^N-\rho^0(x,y)\right|  \leq  \frac{\max(\Delta_1,\Delta_2)}{N} (\|\partial_x \rho^0\|_\infty +
\|\partial_y \rho^0\|_\infty),
$$
then
\begin{align*}
& \left|   \sum_{i=0}^{N-1} \sum_{j=0}^{N-1} \int_{\overline{x}_i^N}^{\overline{x}_{i+1}^N}
\int_{\overline{y}_j^N}^{\overline{y}_{j+1}^N} (a_{i,j}^N-\rho^0(x,y)) ( \varphi(x,y) - \varphi(\overline{x}_i^N,\overline{y}_j^N))\,  \,dydx \right| \\
\leq &  \sum_{i=0}^{N-1} \sum_{j=0}^{N-1} \int_{\overline{x}_i^N}^{\overline{x}_{i+1}^N}
\int_{\overline{y}_j^N}^{\overline{y}_{j+1}^N} |a_{i,j}^N-\rho^0(x,y)| | \varphi(x,y) - \varphi(\overline{x}_i^N,\overline{y}_j^N)| \,  \,dydx \\
\leq &  \frac{\max(\Delta_1,\Delta_2)}{N} (\|\partial_x \rho^0\|_\infty +
\|\partial_y \rho^0\|_\infty) \sum_{i=0}^{N-1} \sum_{j=0}^{N-1} \int_{\overline{x}_i^N}^{\overline{x}_{i+1}^N}
\int_{\overline{y}_j^N}^{\overline{y}_{j+1}^N}  | \varphi(x,y) - \varphi(\overline{x}_i^N,\overline{y}_j^N)| \,  \,dydx \\
\leq &  \frac{\max(\Delta_1,\Delta_2)}{N} (\|\partial_x \rho^0\|_\infty +
\|\partial_y \rho^0\|_\infty)
 ( \|\partial_x \varphi\|_\infty+ \|\partial_y \varphi\|_\infty) \frac{\max(\Delta_1,\Delta_2)^3}{N}
\end{align*}
by using Lemma \ref{lemmaapprox3preli}.
Finally we get the result. \CQFD

We now split the study according to the hypothesis on the support, compact from non compact.
In the case 1, that is to say 
when $\rho^0 \in L^1(\r^2)$ has a compact support included
in some $C_L$ where we chose this $L$ for the location area of
particles.
Then we have $\ds  \iint_{\r^2\setminus C_{L}} \rho^0(x,y) \varphi(x,y) \,dxdy=0$ and
the Proposition \ref{propapprox1} leads to
 the following result.
\begin{Prop} \label{CV_W1i_SC}
Let $\rho^0 \in L^1(\r^2)$ such that $\rho^0 \geq 0$ 
with a compact support included in a $C_{L}$
where $L=(L_1, L_1',L_2,L_2')\in\r^4$ be such that
(\ref{hypotheseL}).
We assume furthermore that $\partial_x \rho^0, \partial_y \rho^0 \in L^\infty(\r^2)$.
For any $N\in \N^*$, we consider
$\hat{\rho}^N$ defined by (\ref{densitediscrete}) with 
(\ref{defpreli1})-(\ref{defpreli2}).
Then for any function $\varphi \in C(\r^2)\cap L^\infty(\r^2)$ such that $\partial_x \varphi, \partial_y \varphi \in L^\infty(\r^2)$,
\begin{equation} \label{Estimation1aS}
  \left| \iint_{\r^2}  \hat{\rho}^{N}(x,y) \varphi(x,y) \,dxdy  - \iint_{\r^2} \rho^0(x,y) \varphi(x,y) \,dxdy \right| 
\leq    \frac{C_{1,2}(\varphi)}{N^2},
\end{equation}
with
\begin{equation} \label{C12S}
C_{1,2}(\varphi) =   \max(\Delta_1,\Delta_2)^4 (\|\partial_x \rho^0\|_\infty +
\|\partial_y \rho^0\|_\infty)
 ( \|\partial_x \varphi\|_\infty+ \|\partial_y \varphi\|_\infty).
\end{equation}
\end{Prop}
This is the first part of Theorem \ref{TH1}.

In the case 2 of the initial data, that is to say where $\rho^0 \in L^1(\r^2)$ doesn't not a compact support,
we fix the location area of the particles by removing a mass 
smaller than an epsilon from this area.
For any $\eps>0$, there exists $L>0$ such that (\ref{hypotheseepsL}).
Then we set $L_1=L_2=L$, $L_1'=L_2'=-L$ which gives $C_{L}=[-L,L]^2$.
Remember also that in this case $\hat{\rho}^{N}$ depends on $\eps$ via $L$
thus we note $\hat{\rho}^{N,\eps}$ instead.
We have then the following result.
\begin{Prop}  \label{CV_W1i_SNC}
Let $\rho^0 \in L^1(\r^2)$ and such that $\rho^0 \geq 0$.
Let $\eps >0$.
Consider $L>0$ such that (\ref{hypotheseepsL}).
We assume furthermore that $\partial_x \rho^0, \partial_y \rho^0 \in L^\infty(\r^2)$.
For any $N\in \N^*$, we consider
$\hat{\rho}^{N,\eps}$ defined by (\ref{densitediscrete}) with 
(\ref{defpreli1})-(\ref{defpreli2}).
Then for any function $\varphi \in C(\r^2)\cap L^\infty(\r^2)$ such that $\partial_x \varphi, \partial_y \varphi \in L^\infty(\r^2)$,
we have
\begin{equation} \label{Estimation2aS}
\left|  \iint_{\r^2} \hat{\rho}^{N,\eps}(x,y) \varphi(x,y) \,dxdy  - \iint_{\r^2} \rho^0(x,y) \varphi(x,y) \,dxdy  \right| \leq   \eps  \|\varphi \|_\infty  + \frac{C_\eps(\varphi)}{N^2},
\end{equation}
with
\begin{equation} \label{CepsS}
C_\eps(\varphi) = 16 L^4 (\|\partial_x \rho^0\|_\infty +
\|\partial_y \rho^0\|_\infty)
 ( \|\partial_x \varphi\|_\infty+ \|\partial_y \varphi\|_\infty).
\end{equation}
\end{Prop}
\pr  
Thanks to Proposition \ref{propapprox1}, we have
\begin{align*}
&  \left| \iint_{\r^2}  \hat{\rho}^{N,\eps}(x,y) \varphi(x,y) \,dxdy  - \iint_{\r^2} \rho^0(x,y) \varphi(x,y) \,dxdy \right| \\
\leq & \left| \iint_{\r^2\setminus C_{L}} \rho^0(x,y) \varphi(x,y) \,dxdy\right|
 +  \frac{\max(\Delta_1,\Delta_2)^4}{N^2} (\|\partial_x \rho^0\|_\infty +
\|\partial_y \rho^0\|_\infty)
 ( \|\partial_x \varphi\|_\infty+ \|\partial_y \varphi\|_\infty).
\end{align*}
Now 
\begin{align*}
\left| \iint_{\r^2\setminus C_L} \rho^0(x,y) \varphi(x,y) \,dxdy \right| & \leq \| \varphi \|_\infty \iint_{\r^2\setminus C_L} \rho^0(x,y)  \,dxdy
\leq  \| \varphi \|_\infty \eps,
\end{align*}
and $\Delta_1=\Delta_2=2L$,
thus we get the result.
Notice that $C_\eps(\varphi)$ depends on $\eps$ since $L$ depends on $\eps$.
\CQFD

This is the first part of Theorem \ref{TH2}.

\subsection{Case where $\rho^0 \in L^{1}(\r^2)\cap L^{\infty}(\r^2)$}

We move now to the estimate of the density approximation in the case of regularity 
$\rho^0 \in L^{1}(\r^2)\cap L^{\infty}(\r^2)$.
In this case, we cannot use the smallest of the term $a_{i,j}^N-\rho^0(x,y)$
and we have to study more precisely.
Instead of having an estimate on $C/N^2$, we will have only a form as $C/N$.
We start with two preliminary estimates.
\begin{Lemma}  \label{lemmaapprox2}
Let $\rho^0 \in L^1(\r^2)$ such that $\rho^0 \geq 0$.
Let $L=(L_1, L_1',L_2,L_2')\in\r^4$ be such that
(\ref{hypotheseL})
and set $\Delta_1 = L_1-L_1'$, $\Delta_2 = L_2-L_2'$ and $C_{L}=[L_1',L_1]\times[L_2',L_2]$.
Let $N\in \N^*$ and
$(\overline{x}_i^N)_{i \in \ent{0}{N-1}}$, $(\overline{y}_j^N)_{j \in \ent{0}{N-1}}$ and 
$(a_{i,j}^N)_{i \in \ent{0}{N-1}, j \in \ent{0}{N-1}}$
defined by 
(\ref{defpreli1}) and (\ref{defpreli2}).
Then for any function $\varphi \in C(\r^2)\cap L^\infty(\r^2)$ such that $\partial_x \varphi, \partial_y \varphi \in L^\infty(\r^2)$,
we have
\begin{align*}
& \left| \sum_{i=0}^{N-1} \sum_{j=0}^{N-1} a_{i,j}^N 
\int_{\overline{x}_i^N}^{\overline{x}_{i+1}^N} \int_{\overline{y}_j^N}^{\overline{y}_{j+1}^N}  
 ( \varphi(x,y) - \varphi(\overline{x}_i^N,\overline{y}_j^N)) \, \,dydx \right| \\
\leq &     ( \|\partial_x \varphi\|_\infty+ \|\partial_y \varphi\|_\infty) \frac{\max(\Delta_1,\Delta_2)}{N} 
 \iint_{\r^2} \rho^0(x,y) \,dxdy.  
\end{align*}
\end{Lemma}
\pr
We have
\begin{align*}
& \left| \sum_{i=0}^{N-1} \sum_{j=0}^{N-1} a_{i,j}^N 
\int_{\overline{x}_i^N}^{\overline{x}_{i+1}^N} \int_{\overline{y}_j^N}^{\overline{y}_{j+1}^N}  
 ( \varphi(x,y) - \varphi(\overline{x}_i^N,\overline{y}_j^N)) \, \,dydx \right|\\
\leq &  \sum_{i=0}^{N-1} \sum_{j=0}^{N-1} |a_{i,j}^N| \, \left|
\int_{\overline{x}_i^N}^{\overline{x}_{i+1}^N} \int_{\overline{y}_j^N}^{\overline{y}_{j+1}^N}  
  ( \varphi(x,y) - \varphi(\overline{x}_i^N,\overline{y}_j^N)) \,  \,dydx \right|\\
\leq &  \sum_{i=0}^{N-1} \sum_{j=0}^{N-1} |a_{i,j}^N| \, \left|
\int_{\overline{x}_i^N}^{\overline{x}_{i+1}^N} \int_{\overline{y}_j^N}^{\overline{y}_{j+1}^N}  
  ( \varphi(x,y) - \varphi(x,\overline{y}_j^N) + \varphi(x,\overline{y}_j^N) -  \varphi(\overline{x}_i^N,\overline{y}_j^N)) \,  \,dydx \right|,
\end{align*}
thus
\begin{align*}
& \left| \sum_{i=0}^{N-1} \sum_{j=0}^{N-1} a_{i,j}^N 
\int_{\overline{x}_i^N}^{\overline{x}_{i+1}^N} \int_{\overline{y}_j^N}^{\overline{y}_{j+1}^N}  
 ( \varphi(x,y) - \varphi(\overline{x}_i^N,\overline{y}_j^N)) \, \,dydx \right|\\
\leq &  \sum_{i=0}^{N-1} \sum_{j=0}^{N-1} |a_{i,j}^N| \, \left|
\int_{\overline{x}_i^N}^{\overline{x}_{i+1}^N} \left( \int_{\overline{y}_j^N}^{\overline{y}_{j+1}^N}  
  ( \varphi(x,y) - \varphi(x,\overline{y}_j^N) ) \,dy  \right) \,  \,dx \right|\\
 &  + \sum_{i=0}^{N-1} \sum_{j=0}^{N-1} |a_{i,j}^N| \, \left|
 \int_{\overline{y}_j^N}^{\overline{y}_{j+1}^N}  \left( \int_{\overline{x}_i^N}^{\overline{x}_{i+1}^N}
  (  \varphi(x,\overline{y}_j^N) -  \varphi(\overline{x}_i^N,\overline{y}_j^N))  \,dx \right)  \,  dy \right| \\
   \leq &  \sum_{i=0}^{N-1} \sum_{j=0}^{N-1} |a_{i,j}^N| \, 
\int_{\overline{x}_i^N}^{\overline{x}_{i+1}^N} \left( \int_{\overline{y}_j^N}^{\overline{y}_{j+1}^N}  
   \|\partial_y \varphi\|_\infty (y-\overline{y}_j^N) \,dy  \right)   \,dx \\
 &  + \sum_{i=0}^{N-1} \sum_{j=0}^{N-1} |a_{i,j}^N| \, 
 \int_{\overline{y}_j^N}^{\overline{y}_{j+1}^N}  \left( \int_{\overline{x}_i^N}^{\overline{x}_{i+1}^N}
  \|\partial_y \varphi\|_\infty (x-\overline{x}_i^N)   \,dx \right)  \,  dy.
\end{align*}
Now
$$
0 \leq  \int_{\overline{y}_j^N}^{\overline{y}_{j+1}^N}  
   \|\partial_y \varphi\|_\infty (y-\overline{y}_j^N) \,dy 
 \leq  \|\partial_y \varphi\|_\infty \left[ \frac12 (y-\overline{y}_j^N )^2 \right]_{\overline{y}_j^N}^{\overline{y}_{j+1}^N}  
 \leq \|\partial_y \varphi\|_\infty  \frac12\frac{\Delta_2^2}{N^2}  
$$
and
$$
0 \leq  \int_{\overline{x}_i^N}^{\overline{x}_{i+1}^N}  
   \|\partial_x \varphi\|_\infty (x-\overline{x}_i^N) \,dx
 \leq  \|\partial_x \varphi\|_\infty \left[ \frac12 (x-\overline{x}_i^N )^2 \right]_{\overline{x}_i^N}^{\overline{x}_{i+1}^N}  
 \leq \|\partial_x \varphi\|_\infty  \frac12\frac{\Delta_1^2}{N^2},  
$$
then we obtain
\begin{align*}
& \left| \sum_{i=0}^{N-1} \sum_{j=0}^{N-1} a_{i,j}^N 
\int_{\overline{x}_i^N}^{\overline{x}_{i+1}^N} \int_{\overline{y}_j^N}^{\overline{y}_{j+1}^N}  
 ( \varphi(x,y) - \varphi(\overline{x}_i^N,\overline{y}_j^N)) \, \,dydx \right|\\
 \leq & \|\partial_y \varphi\|_\infty \sum_{i=0}^{N-1} \sum_{j=0}^{N-1} |a_{i,j}^N| \, 
\frac{\Delta_1}{N}   \frac{\Delta_2^2}{N^2}  
  +  \|\partial_x \varphi\|_\infty \sum_{i=0}^{N-1} \sum_{j=0}^{N-1} |a_{i,j}^N| \, \frac{\Delta_2}{N}   \frac{\Delta_1^2}{N^2}  \\
 \leq &  ( \|\partial_x \varphi\|_\infty+ \|\partial_y \varphi\|_\infty) 
 \frac{\max(\Delta_1,\Delta_2)}{N} 
  \sum_{i=0}^{N-1} \sum_{j=0}^{N-1} 
 \int_{\overline{x}_i^N}^{\overline{x}_{i+1}^N}  \int_{\overline{y}_j^N}^{\overline{y}_{j+1}^N} \rho^0(x,y) \,dydx  \\
\leq   &  ( \|\partial_x \varphi\|_\infty+ \|\partial_y \varphi\|_\infty) \frac{\max(\Delta_1,\Delta_2)}{N} 
 \iint_{\r^2} \rho^0(x,y) \,dxdy.  \CQFD
\end{align*}

\begin{Lemma}  \label{lemmaapprox3}
Let $\rho^0 \in L^1(\r^2)\cap L^{\infty}(\r^2)$ such that $\rho^0 \geq 0$.
Let $L=(L_1, L_1',L_2,L_2')\in\r^4$ be such that
(\ref{hypotheseL})
and set $\Delta_1 = L_1-L_1'$, $\Delta_2 = L_2-L_2'$ and $C_{L}=[L_1',L_1]\times[L_2',L_2]$.
Let $N\in \N^*$ and
$(\overline{x}_i^N)_{i \in \ent{0}{N-1}}$, $(\overline{y}_j^N)_{j \in \ent{0}{N-1}}$
defined by  (\ref{defpreli1}).
Then for any function $\varphi \in C(\r^2)\cap L^\infty(\r^2)$ such that $\partial_x \varphi, \partial_y \varphi \in L^\infty(\r^2)$,
we have
\begin{align*}
& \left| \sum_{i=0}^{N-1} \sum_{j=0}^{N-1} \int_{\overline{x}_i^N}^{\overline{x}_{i+1}^N}
\int_{\overline{y}_j^N}^{\overline{y}_{j+1}^N} \rho^0(x,y) (  \varphi(x,y) - \varphi(\overline{x}_i^N,\overline{y}_j^N)) \,dydx \right|\\
\leq &    \|\rho^0\|_\infty  ( \|\partial_x \varphi\|_\infty+ \|\partial_y \varphi\|_\infty) \frac{\Delta_1 \Delta_2 \max(\Delta_1,\Delta_2)}{N}.  
\end{align*}
\end{Lemma}
\pr We have
\begin{align*}
& \left| \sum_{i=0}^{N-1} \sum_{j=0}^{N-1} \int_{\overline{x}_i^N}^{\overline{x}_{i+1}^N}
\int_{\overline{y}_j^N}^{\overline{y}_{j+1}^N} \rho^0(x,y) (  \varphi(x,y) - \varphi(\overline{x}_i^N,\overline{y}_j^N))  \,dydx \right|\\
\leq &  \|\rho^0\|_\infty  \sum_{i=0}^{N-1} \sum_{j=0}^{N-1}  \, 
\int_{\overline{x}_i^N}^{\overline{x}_{i+1}^N} \int_{\overline{y}_j^N}^{\overline{y}_{j+1}^N}  
  \left| \varphi(x,y) - \varphi(\overline{x}_i^N,\overline{y}_j^N) \right|  \,dydx\\
 \leq & \|\rho^0\|_\infty  ( \|\partial_x \varphi\|_\infty+ \|\partial_y \varphi\|_\infty) \frac{\Delta_1 \Delta_2 \max(\Delta_1,\Delta_2)}{N}.  
\end{align*}
by using Lemma \ref{lemmaapprox3preli}.
\CQFD

We will now consider the two cases of locations for the  particles according to the assumptions on the initial data.
In the case 1, that is to say 
when $\rho^0 \in L^1(\r^2)$ has a compact support included
in some $C_L$ where we choose this $L$ for the location area of
particles, we have the following result.
\begin{Prop} \label{CV_L1_SC}
Let $\rho^0 \in L^1(\r^2)\cap L^{\infty}(\r^2)$ such that $\rho^0 \geq 0$ 
with a compact support included in a $C_{L}$
where $L=(L_1, L_1',L_2,L_2')\in\r^4$ be such that
(\ref{hypotheseL}).
For any $N\in \N^*$, we consider
$\hat{\rho}^N$ defined by (\ref{densitediscrete}) with 
(\ref{defpreli1})-(\ref{defpreli2}).
Then for any function $\varphi \in C(\r^2)\cap L^\infty(\r^2)$ such that $\partial_x \varphi, \partial_y \varphi \in L^\infty(\r^2)$,
we have
\begin{equation} \label{Estimation3aS}
    \left|  \iint_{\r^2} \hat{\rho}^N(x,y) \varphi(x,y) \,dxdy  - \iint_{\r^2} \rho^0(x,y) \varphi(x,y) \,dxdy  \right| \leq  \frac{D_{1,2}(\varphi)}{N}.
\end{equation}
with
\begin{equation} \label{D12S}
{D}_{1,2}(\varphi) =  \max(\Delta_1,\Delta_2) 
 ( \|\partial_x \varphi\|_\infty+ \|\partial_y \varphi\|_\infty)
 (\|\rho^0\|_1 + \|\rho^0\|_\infty \Delta_1\Delta_2).
\end{equation}
\end{Prop}
\pr 
Thanks to Lemma \ref{lemmaapprox1}, we have
\begin{align*}
&   \iint_{\r^2} \hat{\rho}^N(x,y) \varphi(x,y) \,dxdy  - \iint_{\r^2} \rho^0(x,y) \varphi(x,y) \,dxdy  \\
= & - \iint_{\r^2\setminus C_{L}} \rho^0(x,y) \varphi(x,y) \,dxdy +  \sum_{i=0}^{N-1} \sum_{j=0}^{N-1} \int_{\overline{x}_i^N}^{\overline{x}_{i+1}^N}
\int_{\overline{y}_j^N}^{\overline{y}_{j+1}^N} (a_{i,j}^N-\rho^0(x,y)) ( \varphi(x,y) - \varphi(\overline{x}_i^N,\overline{y}_j^N))\,  \,dydx.
\end{align*}
Since $\rho^0$ have a compact support in $C_{L}$, we get
$\ds \iint_{\r^2\setminus C_{L}} \rho^0(x,y) \varphi(x,y) \,dxdy=0$.
With Lemma \ref{lemmaapprox2}, we have
\begin{align*}
& \left| \sum_{i=0}^{N-1} \sum_{j=0}^{N-1} a_{i,j}^N 
\int_{\overline{x}_i^N}^{\overline{x}_{i+1}^N} \int_{\overline{y}_j^N}^{\overline{y}_{j+1}^N}  
 ( \varphi(x,y) - \varphi(\overline{x}_i^N,\overline{y}_j^N)) \, \,dydx \right|\\
\leq   &  ( \|\partial_x \varphi\|_\infty+ \|\partial_y \varphi\|_\infty) \frac{\max(\Delta_1,\Delta_2)}{N} 
 \iint_{\r^2} \rho^0(x,y) \,dxdy.  
\end{align*}
and with Lemma \ref{lemmaapprox3}, we have
$$
 \left| \sum_{i=0}^{N-1} \sum_{j=0}^{N-1} \int_{\overline{x}_i^N}^{\overline{x}_{i+1}^N}
\int_{\overline{y}_j^N}^{\overline{y}_{j+1}^N} \rho^0(x,y) (  \varphi(x,y) - \varphi(\overline{x}_i^N,\overline{y}_j^N)) \,dydx \right|
\leq     \|\rho^0\|_\infty  ( \|\partial_x \varphi\|_\infty+ \|\partial_y \varphi\|_\infty) \frac{\Delta_1 \Delta_2 \max(\Delta_1,\Delta_2)}{N}.  
$$
Thus we get
\begin{align*}
&  \left| \iint_{\r^2} \hat{\rho}^N(x,y) \varphi(x,y) \,dxdy  - \iint_{\r^2} \rho^0(x,y) \varphi(x,y) \,dxdy \right| \\
\leq & 
 ( \|\partial_x \varphi\|_\infty+ \|\partial_y \varphi\|_\infty) \frac{\max(\Delta_1,\Delta_2)}{N} 
 \iint_{\r^2} \rho^0(x,y) \,dxdy
 +  \|\rho^0\|_\infty  ( \|\partial_x \varphi\|_\infty+ \|\partial_y \varphi\|_\infty) \frac{\Delta_1 \Delta_2 \max(\Delta_1,\Delta_2)}{N}
\end{align*}
and the expected estimate.
\CQFD

This is the first part of Theorem \ref{TH3}.

In the case 2 of the initial data, that is to say where $\rho^0 \in L^1(\r^2)$ doesn't not a compact support.
we fix the location area of the particles by removing a mass 
smaller than an epsilon from this area.
For any $\eps>0$, there exists $L>0$ such that (\ref{hypotheseepsL}).
Then we set $L_1=L_2=L$, $L_1'=L_2'=-L$ which gives $C_{L}=[-L,L]^2$.
Remember also that in this case $\hat{\rho}^{N}$ depends on $\eps$ via $L$
thus we note $\hat{\rho}^{N,\eps}$ instead.
We have the following result.
\begin{Prop}  \label{CV_L1_SNC}
Let $\rho^0 \in L^1(\r^2)\cap L^{\infty}(\r^2)$ and such that $\rho^0 \geq 0$.
Let $\eps >0$.
Consider $L>0$ such that (\ref{hypotheseepsL}).
For any $N\in \N^*$, we consider
$\hat{\rho}^{N,\eps}$ defined by (\ref{densitediscrete}) with 
(\ref{defpreli1})-(\ref{defpreli2}).
Then for any function $\varphi \in C(\r^2)\cap L^\infty(\r^2)$ such that $\partial_x \varphi, \partial_y \varphi \in L^\infty(\r^2)$,
we have
\begin{equation} \label{Estimation4aS}
\left|  \iint_{\r^2} \hat{\rho}^{N,\eps}(x,y) \varphi(x,y) \,dxdy  - \iint_{\r^2} \rho^0(x,y) \varphi(x,y) \,dxdy  \right| \leq \eps \|  \varphi \|_\infty  + \frac{D_\eps(\varphi)}{N},
\end{equation}
with
\begin{equation} \label{DepsS}
D_\eps(\varphi) = 2L ( \|\partial_x \varphi\|_\infty+ \|\partial_y \varphi\|_\infty) 
( \|\rho^0\|_1  +    4L^2  \|\rho^0\|_\infty).
\end{equation}
\end{Prop}
\pr  
Thanks to Lemma \ref{lemmaapprox1}, we have
\begin{align*}
&   \iint_{\r^2}  \hat{\rho}^{N,\eps}(x,y) \varphi(x,y) \,dxdy  - \iint_{\r^2} \rho^0(x,y) \varphi(x,y) \,dxdy  \\
= & - \iint_{\r^2\setminus C_{L}} \rho^0(x,y) \varphi(x,y) \,dxdy +  \sum_{i=0}^{N-1} \sum_{j=0}^{N-1} \int_{\overline{x}_i^N}^{\overline{x}_{i+1}^N}
\int_{\overline{y}_j^N}^{\overline{y}_{j+1}^N} (a_{i,j}^N-\rho^0(x,y)) ( \varphi(x,y) - \varphi(\overline{x}_i^N,\overline{y}_j^N))\,  \,dydx.
\end{align*}
First we have
\begin{align*}
\left| \iint_{\r^2\setminus C_L} \rho^0(x,y) \varphi(x,y) \,dxdy \right| & \leq \| \varphi \|_\infty \iint_{\r^2\setminus C_L} \rho^0(x,y)  \,dxdy
\leq  \| \varphi \|_\infty \eps.
\end{align*}
With Lemma \ref{lemmaapprox2}, we have
\begin{align*}
& \left| \sum_{i=0}^{N-1} \sum_{j=0}^{N-1} a_{i,j}^N 
\int_{\overline{x}_i^N}^{\overline{x}_{i+1}^N} \int_{\overline{y}_j^N}^{\overline{y}_{j+1}^N}  
 ( \varphi(x,y) - \varphi(\overline{x}_i^N,\overline{y}_j^N)) \, \,dydx \right|\\
\leq   &  ( \|\partial_x \varphi\|_\infty+ \|\partial_y \varphi\|_\infty) \frac{\Delta}{N} 
 \iint_{\r^2} \rho^0(x,y) \,dxdy.  
\end{align*}
where $\Delta=2L$,
and with Lemma \ref{lemmaapprox3}, we have
$$
 \left| \sum_{i=0}^{N-1} \sum_{j=0}^{N-1} \int_{\overline{x}_i^N}^{\overline{x}_{i+1}^N}
\int_{\overline{y}_j^N}^{\overline{y}_{j+1}^N} \rho^0(x,y) (  \varphi(x,y) - \varphi(\overline{x}_i^N,\overline{y}_j^N)) \,dydx \right|
\leq     \|\rho^0\|_\infty  ( \|\partial_x \varphi\|_\infty+ \|\partial_y \varphi\|_\infty) \frac{\Delta^3}{N}.  
$$
Thus we get
\begin{align*}
&  \left| \iint_{\r^2}  \hat{\rho}^{N,\eps}(x,y) \varphi(x,y) \,dxdy  - \iint_{\r^2} \rho^0(x,y) \varphi(x,y) \,dxdy \right| \\
\leq & \| \varphi \|_\infty \eps + 
 ( \|\partial_x \varphi\|_\infty+ \|\partial_y \varphi\|_\infty) \frac{\Delta}{N} 
 \iint_{\r^2} \rho^0(x,y) \,dxdy
 +  \|\rho^0\|_\infty  ( \|\partial_x \varphi\|_\infty+ \|\partial_y \varphi\|_\infty) \frac{\Delta^3}{N}
\end{align*}
and finally the wanted inequality since $\Delta=2L$.
Notice that $D_\eps(\varphi)$ depends on $\eps$ since $L$ depends on $\eps$.
\CQFD

This is the first part of Theorem \ref{TH4}.

\section{Approximation of other quantities}

\subsection{General estimates}

We turn now to the estimate of additionnal quantities of shape 
$\ds \iint_{\r^2} \rho^0(x,y) \omega^0(x,y) \varphi(x,y) \,dxdy$.
We start by presenting properties valid in every cases of regularity and support on density $\rho^0$.
\begin{Lemma} \label{lemmaapproxG1}
Let $\rho^0 \in L^1(\r^2)$ such that $\rho^0 \geq 0$ and let $\omega^0 \in L^\infty(\r^2)$.
Let $L=(L_1, L_1',L_2,L_2')\in\r^4$ be such that
(\ref{hypotheseL})
and set $\Delta_1 = L_1-L_1'$, $\Delta_2 = L_2-L_2'$ and $C_{L}=[L_1',L_1]\times[L_2',L_2]$.
Let $N\in \N^*$ and
$\hat{\rho}^N, \hat{\rho}^N \hat{\omega}^N$ defined by (\ref{densitediscrete}) and (\ref{quantitydiscrete}) with 
(\ref{defpreli1})-(\ref{defpreli2}) and (\ref{initialdiscretquantity}).
Then for any function $\varphi \in C(\r^2)\cap L^\infty(\r^2)$,
we have
\begin{align*}
&   \iint_{\r^2} \hat{\rho}^N(x,y) \hat{\omega}^N(x,y) \varphi(x,y) \,dxdy  - \iint_{\r^2} \rho^0(x,y)
\omega^0(x,y) \varphi(x,y) \,dxdy  \\
= & - \iint_{\r^2\setminus C_{L}} \rho^0(x,y) \omega^0(x,y) 
\varphi(x,y) \,dxdy +  \sum_{i=0}^{N-1} \sum_{j=0}^{N-1} \int_{\overline{x}_i^N}^{\overline{x}_{i+1}^N}
\int_{\overline{y}_j^N}^{\overline{y}_{j+1}^N} \rho^0(x,y) (W_{i,j}^N -\omega^0(x,y)) \varphi(\overline{x}_i^N,\overline{y}_j^N)\,  \,dydx\\
 & +  \sum_{i=0}^{N-1} \sum_{j=0}^{N-1} W_{i,j}^N \int_{\overline{x}_i^N}^{\overline{x}_{i+1}^N}
\int_{\overline{y}_j^N}^{\overline{y}_{j+1}^N} (a_{i,j}^N-\rho^0(x,y)) ( \varphi(x,y) - \varphi(\overline{x}_i^N,\overline{y}_j^N))\,  \,dydx.
\end{align*}
\end{Lemma}
\pr
Since $\varphi, \omega^0 \in L^\infty(\r^2)$ and
$\rho^0 \in L^1(\r^2)$, then
$\ds \iint_{\r^2} \rho^0(x,y) \omega^0(x,y) \varphi(x,y) \,dxdy$ exists.
Notice that with the definition of $\hat{\rho}^N$, we have $\hat{\rho}^N=0$ on $\r^2\setminus C_{L}$.
Doing similarly to the proof of Lemma \ref{lemmaapprox1}, we obtain
\begin{align*}
 \iint_{\r^2} \hat{\rho}^N(x,y) \hat{\omega}^N(x,y) \varphi(x,y) \,dxdy  
= &  
\sum_{i=0}^{N-1} \sum_{j=0}^{N-1} 
\int_{\overline{x}_i^N}^{\overline{x}_{i+1}^N} \int_{\overline{y}_j^N}^{\overline{y}_{j+1}^N} 
 a_{i,j}^N W_{i,j}^N \varphi(\overline{x}_i^N,\overline{y}_j^N) \,  \,dydx\\
&  +\sum_{i=0}^{N-1} \sum_{j=0}^{N-1} a_{i,j}^N W_{i,j}^N
\int_{\overline{x}_i^N}^{\overline{x}_{i+1}^N} \int_{\overline{y}_j^N}^{\overline{y}_{j+1}^N}  
 ( \varphi(x,y) - \varphi(\overline{x}_i^N,\overline{y}_j^N)) \,  \,dydx,
\end{align*}
and
\begin{align*}
& \iint_{\r^2} \rho^0(x,y) \omega^0(x,y) \varphi(x,y) \,dxdy  \\
= &    \iint_{\r^2\setminus C_{L}} \rho^0(x,y)  \omega^0(x,y) \varphi(x,y) \,dxdy + 
\sum_{i=0}^{N-1} \sum_{j=0}^{N-1} 
\int_{\overline{x}_i^N}^{\overline{x}_{i+1}^N} \int_{\overline{y}_j^N}^{\overline{y}_{j+1}^N} \rho^0(x,y) 
 \omega^0(x,y)  \varphi(\overline{x}_i^N,\overline{y}_j^N)\,  \,dydx \\
&  + \sum_{i=0}^{N-1} \sum_{j=0}^{N-1} \int_{\overline{x}_i^N}^{\overline{x}_{i+1}^N}
\int_{\overline{y}_j^N}^{\overline{y}_{j+1}^N} \rho^0(x,y)  \omega^0(x,y) 
(  \varphi(x,y) - \varphi(\overline{x}_i^N,\overline{y}_j^N))\,  \,dydx
\end{align*}
and also
\begin{align*}
 \sum_{i=0}^{N-1} \sum_{j=0}^{N-1} 
\int_{\overline{x}_i^N}^{\overline{x}_{i+1}^N} \int_{\overline{y}_j^N}^{\overline{y}_{j+1}^N}  a_{i,j}^N  
 W_{i,j}^N \varphi(\overline{x}_i^N,\overline{y}_j^N) \,  \,dydx 
= &   \sum_{i=0}^{N-1} \sum_{j=0}^{N-1}  W_{i,j}^N
  \int_{\overline{x}_i^N}^{\overline{x}_{i+1}^N}  \int_{\overline{y}_j^N}^{\overline{y}_{j+1}^N}  \rho^0(x,y)  \varphi(\overline{x}_i^N,\overline{y}_j^N)
  \,  \,dydx.
\end{align*}
Thus we get
\begin{align*}
&   \iint_{\r^2} \hat{\rho}^N(x,y) \hat{\omega}^N(x,y)  \varphi(x,y) \,dxdy 
 - \iint_{\r^2} \rho^0(x,y)  \omega^0(x,y) \varphi(x,y) \,dxdy  \\
= & - \iint_{\r^2\setminus C_{L}} \rho^0(x,y)  \omega^0(x,y) \varphi(x,y) \,dxdy +
\sum_{i=0}^{N-1} \sum_{j=0}^{N-1}  W_{i,j}^N
  \int_{\overline{x}_i^N}^{\overline{x}_{i+1}^N}  \int_{\overline{y}_j^N}^{\overline{y}_{j+1}^N}  \rho^0(x,y)  \varphi(\overline{x}_i^N,\overline{y}_j^N)
  \,  \,dydx \\
& + \sum_{i=0}^{N-1} \sum_{j=0}^{N-1} a_{i,j}^N W_{i,j}^N
\int_{\overline{x}_i^N}^{\overline{x}_{i+1}^N} \int_{\overline{y}_j^N}^{\overline{y}_{j+1}^N}  
 ( \varphi(x,y) - \varphi(\overline{x}_i^N,\overline{y}_j^N)) \,  \,dydx \\
 & - \sum_{i=0}^{N-1} \sum_{j=0}^{N-1} 
\int_{\overline{x}_i^N}^{\overline{x}_{i+1}^N} \int_{\overline{y}_j^N}^{\overline{y}_{j+1}^N} \rho^0(x,y) 
 \omega^0(x,y)  \varphi(\overline{x}_i^N,\overline{y}_j^N)\,  \,dydx \\
&  - \sum_{i=0}^{N-1} \sum_{j=0}^{N-1} \int_{\overline{x}_i^N}^{\overline{x}_{i+1}^N}
\int_{\overline{y}_j^N}^{\overline{y}_{j+1}^N} \rho^0(x,y)  \omega^0(x,y) 
(  \varphi(x,y) - \varphi(\overline{x}_i^N,\overline{y}_j^N))\,  \,dydx
\end{align*}
and the result.
\CQFD

\begin{Lemma}  \label{lemmaapproxsupplG1}
Let $\rho^0 \in L^1(\r^2)$ such that $\rho^0 \geq 0$ and let $\omega^0 \in L^\infty(\r^2)$.
Let $L=(L_1, L_1',L_2,L_2')\in\r^4$ be such that
(\ref{hypotheseL})
and set $\Delta_1 = L_1-L_1'$, $\Delta_2 = L_2-L_2'$ and $C_{L}=[L_1',L_1]\times[L_2',L_2]$.
We assume furthermore that $\partial_x \omega^0, \partial_y \omega^0 \in L^\infty(\r^2)$.
Let $N\in \N^*$ and
$(\overline{x}_i^N)_{i \in \ent{0}{N-1}}$, $(\overline{y}_j^N)_{j \in \ent{0}{N-1}}$ and 
$(a_{i,j}^N)_{i \in \ent{0}{N-1}, j \in \ent{0}{N-1}}$
defined by 
(\ref{defpreli1}) and (\ref{defpreli2}).
Then for any function $\varphi \in C(\r^2)\cap L^\infty(\r^2)$ such that $\partial_x \varphi, \partial_y \varphi \in L^\infty(\r^2)$,
we have
\begin{align*}
& \left| \sum_{i=0}^{N-1} \sum_{j=0}^{N-1} \int_{\overline{x}_i^N}^{\overline{x}_{i+1}^N}
\int_{\overline{y}_j^N}^{\overline{y}_{j+1}^N} \rho^0(x,y) (W_{i,j}^N -\omega^0(x,y)) \varphi(\overline{x}_i^N,\overline{y}_j^N)\,  \,dydx \right| \\
\leq &  \frac{\max(\Delta_1,\Delta_2)}{N} (\|\partial_x \omega^0\|_\infty +
\|\partial_y \omega^0\|_\infty) \|\varphi\|_\infty \iint_{\r^2} \rho^0(x,y) \,dxdy.
\end{align*}
\end{Lemma}
\pr 
Similarly to Lemma \ref{lemmaapprox1suppl}, we have
\begin{align*}
\left|W_{i,j}^N-\omega^0(x,y)\right|  \leq  \frac{\max(\Delta_1,\Delta_2)}{N} (\|\partial_x \omega^0\|_\infty +
\|\partial_y \omega^0\|_\infty).
\end{align*}
Then we get
\begin{align*}
& \left| \sum_{i=0}^{N-1} \sum_{j=0}^{N-1} \int_{\overline{x}_i^N}^{\overline{x}_{i+1}^N}
\int_{\overline{y}_j^N}^{\overline{y}_{j+1}^N} \rho^0(x,y) (W_{i,j}^N -\omega^0(x,y)) \varphi(\overline{x}_i^N,\overline{y}_j^N)\,  \,dydx \right| \\
\leq &  \frac{\max(\Delta_1,\Delta_2)}{N} (\|\partial_x \omega^0\|_\infty +
\|\partial_y \omega^0\|_\infty) \|\varphi\|_\infty \sum_{i=0}^{N-1} \sum_{j=0}^{N-1} \int_{\overline{x}_i^N}^{\overline{x}_{i+1}^N}
\int_{\overline{y}_j^N}^{\overline{y}_{j+1}^N} \rho^0(x,y) \,dxdy
\end{align*}
and the result.
\CQFD

\subsection{Case where $\rho^0 \in W^{1,\infty}(\r^2)$}

We now move on to the estimate of the density approximation in the case of regularity 
$\rho^0 \in W^{1,\infty}(\r^2)$.
In this case, the estimate of the difference between $a_{i,j}^N$
and $\rho^0(x,y)$ for $(x,y)\in [\overline{x}_i^N, \overline{x}_{i+1}^N]\times [\overline{y}_j^N,\overline{y}_{j+1}^N]$ allows to get a precise estimate.
We first have the result.
\begin{Prop} \label{propapproxG1}
Let $\rho^0 \in L^1(\r^2)$ such that $\rho^0 \geq 0$ and let $\omega^0 \in L^\infty(\r^2)$.
We assume furthermore that $\partial_x \rho^0, \partial_y \rho^0, 
\partial_x \omega^0, \partial_y \omega^0 \in L^\infty(\r^2)$.
Let $L=(L_1, L_1',L_2,L_2')\in\r^4$ be such that
(\ref{hypotheseL})
and set $\Delta_1 = L_1-L_1'$, $\Delta_2 = L_2-L_2'$ and $C_{L}=[L_1',L_1]\times[L_2',L_2]$.
Let $N\in \N^*$ and
$\hat{\rho}^N, \hat{\rho}^N \hat{\omega}^N$ defined by (\ref{densitediscrete}) and (\ref{quantitydiscrete}) with 
(\ref{defpreli1})-(\ref{defpreli2}) and (\ref{initialdiscretquantity}).
Then for any function $\varphi \in C(\r^2)\cap L^\infty(\r^2)$ such that $\partial_x \varphi, \partial_y \varphi \in L^\infty(\r^2)$,
we have
\begin{align*}
&  \left|  \iint_{\r^2} \hat{\rho}^N(x,y) \hat{\omega}^N(x,y) \varphi(x,y) \,dxdy  - \iint_{\r^2} \rho^0(x,y)
\omega^0(x,y) \varphi(x,y) \,dxdy  \right| \\
\leq & \left| \iint_{\r^2\setminus C_{L}} \rho^0(x,y) \omega^0(x,y) \varphi(x,y) \,dxdy\right|
 +   \frac{\max(\Delta_1,\Delta_2)}{N} (\|\partial_x \omega^0\|_\infty +
\|\partial_y \omega^0\|_\infty) \|\varphi\|_\infty \iint_{\r^2} \rho^0(x,y) \,dxdy \\
& +  (\|\partial_x \rho^0\|_\infty +
\|\partial_y \rho^0\|_\infty)
 ( \|\partial_x \varphi\|_\infty+ \|\partial_y \varphi\|_\infty) \frac{\max(\Delta_1,\Delta_2)^4}{N^2}.
\end{align*}
\end{Prop}
\pr 
Thanks to Lemma \ref{lemmaapproxG1}, we have
\begin{align*}
&   \iint_{\r^2} \hat{\rho}^N(x,y) \hat{\omega}^N(x,y) \varphi(x,y) \,dxdy  - \iint_{\r^2} \rho^0(x,y)
\omega^0(x,y) \varphi(x,y) \,dxdy  \\
= & - \iint_{\r^2\setminus C_{L}} \rho^0(x,y) \omega^0(x,y) 
\varphi(x,y) \,dxdy +  \sum_{i=0}^{N-1} \sum_{j=0}^{N-1} \int_{\overline{x}_i^N}^{\overline{x}_{i+1}^N}
\int_{\overline{y}_j^N}^{\overline{y}_{j+1}^N} \rho^0(x,y) (W_{i,j}^N -\omega^0(x,y)) \varphi(\overline{x}_i^N,\overline{y}_j^N)\,  \,dydx\\
 & +  \sum_{i=0}^{N-1} \sum_{j=0}^{N-1} W_{i,j}^N \int_{\overline{x}_i^N}^{\overline{x}_{i+1}^N}
\int_{\overline{y}_j^N}^{\overline{y}_{j+1}^N} (a_{i,j}^N-\rho^0(x,y)) ( \varphi(x,y) - \varphi(\overline{x}_i^N,\overline{y}_j^N))\,  \,dydx.
\end{align*}
On one side, using Lemma \ref{lemmaapprox1suppl}, we have
$$
\left|a_{i,j}^N-\rho^0(x,y)\right|  \leq  \frac{\max(\Delta_1,\Delta_2)}{N} (\|\partial_x \rho^0\|_\infty +
\|\partial_y \rho^0\|_\infty),
$$
and on the other side, using Lemma \ref{lemmaapprox3preli}
we have
$$
  \sum_{i=0}^{N-1} \sum_{j=0}^{N-1} \,
\int_{\overline{x}_i^N}^{\overline{x}_{i+1}^N} \int_{\overline{y}_j^N}^{\overline{y}_{j+1}^N}  
  \left| \varphi(x,y) - \varphi(\overline{x}_i^N,\overline{y}_j^N) \right|  \,dydx
\leq      ( \|\partial_x \varphi\|_\infty+ \|\partial_y \varphi\|_\infty) 
\frac{\Delta_1 \Delta_2 \max(\Delta_1,\Delta_2)}{N}.  
$$
Then putting this together we get
\begin{align*}
& \left|   \sum_{i=0}^{N-1} \sum_{j=0}^{N-1}  W_{i,j}^N \int_{\overline{x}_i^N}^{\overline{x}_{i+1}^N}
\int_{\overline{y}_j^N}^{\overline{y}_{j+1}^N} (a_{i,j}^N-\rho^0(x,y)) ( \varphi(x,y) - \varphi(\overline{x}_i^N,\overline{y}_j^N))\,  \,dydx \right| \\
\leq & \|\omega^0\|_\infty  \sum_{i=0}^{N-1} \sum_{j=0}^{N-1} \int_{\overline{x}_i^N}^{\overline{x}_{i+1}^N}
\int_{\overline{y}_j^N}^{\overline{y}_{j+1}^N} |a_{i,j}^N-\rho^0(x,y)| | \varphi(x,y) - \varphi(\overline{x}_i^N,\overline{y}_j^N)| \,  \,dydx \\
\leq &  \frac{\max(\Delta_1,\Delta_2)}{N} (\|\partial_x \rho^0\|_\infty +
\|\partial_y \rho^0\|_\infty) \sum_{i=0}^{N-1} \sum_{j=0}^{N-1} \int_{\overline{x}_i^N}^{\overline{x}_{i+1}^N}
\int_{\overline{y}_j^N}^{\overline{y}_{j+1}^N}  | \varphi(x,y) - \varphi(\overline{x}_i^N,\overline{y}_j^N)| \,  \,dydx \\
\leq &  \frac{\max(\Delta_1,\Delta_2)}{N} (\|\partial_x \rho^0\|_\infty +
\|\partial_y \rho^0\|_\infty)
 ( \|\partial_x \varphi\|_\infty+ \|\partial_y \varphi\|_\infty) \frac{\max(\Delta_1,\Delta_2)^3}{N}.
\end{align*}
Now using Lemma \ref{lemmaapproxsupplG1}, we have
\begin{align*}
& \left| \sum_{i=0}^{N-1} \sum_{j=0}^{N-1} \int_{\overline{x}_i^N}^{\overline{x}_{i+1}^N}
\int_{\overline{y}_j^N}^{\overline{y}_{j+1}^N} \rho^0(x,y) (W_{i,j}^N -\omega^0(x,y)) \varphi(\overline{x}_i^N,\overline{y}_j^N)\,  \,dydx \right| \\
\leq &  \frac{\max(\Delta_1,\Delta_2)}{N} (\|\partial_x \omega^0\|_\infty +
\|\partial_y \omega^0\|_\infty) \|\varphi\|_\infty \iint_{\r^2} \rho^0(x,y) \,dxdy
\end{align*}
and we get the result.
\CQFD

We will now consider both options of locations for the  particles according to the assumptions on the initial data.
In the case 1, that is to say 
where $\rho^0 \in L^1(\r^2)$ has a compact support included
in some $C_L$, Proposition \ref{propapproxG1} gives
 the following result.
\begin{Prop} \label{CV_W1i_SC_G}
Let $\rho^0 \in L^1(\r^2)$ such that $\rho^0 \geq 0$
with a compact support included in a $C_{L}$
where $L=(L_1, L_1',L_2,L_2')\in\r^4$ be such that
(\ref{hypotheseL}) and let $\omega^0 \in L^\infty(\r^2)$.
We assume furthermore that $\partial_x \rho^0, \partial_y \rho^0, \partial_x \omega^0, 
\partial_y \omega^0 \in L^\infty(\r^2)$.
Let $N\in \N^*$ and
$\hat{\rho}^N, \hat{\rho}^N \hat{\omega}^N$ defined by (\ref{densitediscrete}) and (\ref{quantitydiscrete}) with 
(\ref{defpreli1})-(\ref{defpreli2}) and (\ref{initialdiscretquantity}).
Then for any function $\varphi \in C(\r^2)\cap L^\infty(\r^2)$ such that $\partial_x \varphi, \partial_y \varphi \in L^\infty(\r^2)$,
we have
\begin{equation} \label{Estimation1bS}
  \left|  \iint_{\r^2} \hat{\rho}^N(x,y) \hat{\omega}^N(x,y) \varphi(x,y) \,dxdy  - \iint_{\r^2} \rho^0(x,y)
\omega^0(x,y) \varphi(x,y) \,dxdy  \right| 
\leq     \frac{K_{1,2}(\varphi)}{N} +\frac{C_{1,2}(\varphi) \|\omega^0\|_\infty}{N^2},
\end{equation}
with 
\begin{equation} \label{K12bS}
K_{1,2}(\varphi) =
\max(\Delta_1,\Delta_2) (\|\partial_x \omega^0\|_\infty +
\|\partial_y \omega^0\|_\infty) \|\varphi\|_\infty \| \rho^0 \|_1
\end{equation}
and $C_{1,2}(\varphi)$ defined by (\ref{C12S}).
\end{Prop}
This proves the remaining part of Theorem \ref{TH1}.
In the case 2 of the initial data, that is to say where $\rho^0 \in L^1(\r^2)$ doesn't not a compact support,
We have the following result.
\begin{Prop}  \label{CV_W1i_SNC_G}
Let $\rho^0 \in L^1(\r^2)$ such that $\rho^0 \geq 0$ and let $\omega^0 \in L^\infty(\r^2)$.
Let $\eps >0$.
Consider $L>0$ such that (\ref{hypotheseepsL}).
We assume furthermore that $\partial_x \rho^0, \partial_y \rho^0, \partial_x \omega^0, 
\partial_y \omega^0 \in L^\infty(\r^2)$.
Let $N\in \N^*$ and
$\hat{\rho}^{N,\eps}, \hat{\rho}^{N,\eps} \hat{\omega}^{N,\eps}$ defined by (\ref{densitediscrete}) and (\ref{quantitydiscrete}) with 
(\ref{defpreli1})-(\ref{defpreli2}) and (\ref{initialdiscretquantity}).
Then for any function $\varphi \in C(\r^2)\cap L^\infty(\r^2)$ such that $\partial_x \varphi, \partial_y \varphi \in L^\infty(\r^2)$,
we have
\begin{align} \nonumber
& \left| \iint_{\r^2} \hat{\rho}^{N,\eps}(x,y) \hat{\omega}^{N,\eps}(x,y) \varphi(x,y) \,dxdy  -
 \iint_{\r^2} \rho^0(x,y) \omega^0(x,y) \varphi(x,y) \,dxdy \right| \\
& \qquad\qquad\qquad \leq   \eps  \|\varphi \|_\infty \|\omega^0\|_\infty + \frac{C_\eps(\varphi)\|\omega^0\|_\infty}{N^2} + \frac{K_\eps(\varphi)}{N},  \label{Estimation2bS}
\end{align}
with 
\begin{equation} \label{KepsS}
K_\eps(\varphi) = 2 L (\|\partial_x \omega^0\|_\infty +
\|\partial_y \omega^0\|_\infty) \|\varphi\|_\infty \|\rho^0\|_1,
\end{equation}
and $C_\eps(\varphi)$ defined by (\ref{CepsS}).
\end{Prop}
\pr  
Thanks to Proposition \ref{propapproxG1}, we have
\begin{align*}
&  \left|  \iint_{\r^2} \hat{\rho}^{N,\eps}(x,y) \hat{\omega}^{N,\eps}(x,y) \varphi(x,y) \,dxdy  - \iint_{\r^2} \rho^0(x,y)
\omega^0(x,y) \varphi(x,y) \,dxdy  \right| \\
\leq & \left| \iint_{\r^2\setminus C_{L}} \rho^0(x,y) \omega^0(x,y) \varphi(x,y) \,dxdy\right|
 +   \frac{\max(\Delta_1,\Delta_2)}{N} (\|\partial_x \omega^0\|_\infty +
\|\partial_y \omega^0\|_\infty) \|\varphi\|_\infty \iint_{\r^2} \rho^0(x,y) \,dxdy \\
& +  (\|\partial_x \rho^0\|_\infty +
\|\partial_y \rho^0\|_\infty)
 ( \|\partial_x \varphi\|_\infty+ \|\partial_y \varphi\|_\infty) \frac{\max(\Delta_1,\Delta_2)^4}{N^2}.
\end{align*}
Now 
\begin{align*}
\left| \iint_{\r^2\setminus C_L} \rho^0(x,y) \omega^0(x,y) \varphi(x,y) \,dxdy \right| & \leq
 \| \varphi \|_\infty \| \omega \|_\infty \iint_{\r^2\setminus C_L} \rho^0(x,y)  \,dxdy
\leq  \| \varphi \|_\infty \| \omega \|_\infty \eps,
\end{align*}
and $\Delta_1=\Delta_2=2L$,
thus we get the result.
Notice that $D_\eps(\varphi)$ depends on $\eps$ since $L$ depends on $\eps$.
\CQFD

This proves the remaining part of Theorem \ref{TH2}.

\subsection{Case where $\rho^0 \in L^{1}(\r^2)\cap L^{\infty}(\r^2)$}

Let's now move to the estimate of additionnal quantities 
in the case of regularity $\rho^0 \in L^{1}(\r^2)\cap L^{\infty}(\r^2)$.
In this case, we cannot use the smallest of the terms $a_{i,j}^N-\rho^0(x,y)$
so we have to analyse deeper.
We split the study again according to the support property of density.
In the case 1, that is to say 
where $\rho^0 \in L^1(\r^2)$ has a compact support included
in some $C_L$, we have the following result.
\begin{Prop} \label{CV_L1_SC_G}
Let $\rho^0 \in L^1(\r^2)\cap L^{\infty}(\r^2)$ such that $\rho^0 \geq 0$ 
with a compact support included in a $C_{L}$
where $L=(L_1, L_1',L_2,L_2')\in\r^4$ be such that
(\ref{hypotheseL}) and let $\omega^0 \in L^\infty(\r^2)$.
We assume furthermore that $\partial_x \omega^0, \partial_y \omega^0 \in L^\infty(\r^2)$.
Let $N\in \N^*$ and
$\hat{\rho}^N, \hat{\rho}^N \hat{\omega}^N$ defined by (\ref{densitediscrete}) and (\ref{quantitydiscrete}) with 
(\ref{defpreli1})-(\ref{defpreli2}) and (\ref{initialdiscretquantity}).
Then for any function $\varphi \in C(\r^2)\cap L^\infty(\r^2)$ such that $\partial_x \varphi, \partial_y \varphi \in L^\infty(\r^2)$,
we have
\begin{equation} \label{Estimation3bS}
  \left|  \iint_{\r^2} \hat{\rho}^N(x,y) \hat{\omega}^N(x,y) \varphi(x,y) \,dxdy  - \iint_{\r^2} \rho^0(x,y)
\omega^0(x,y) \varphi(x,y) \,dxdy  \right| 
\leq     \frac{K_{1,2}(\varphi)+D_{1,2}(\varphi) \|\omega^0\|_\infty}{N},
\end{equation}
where $K_{1,2}(\varphi)$ is defined by (\ref{K12bS}) anf $D_{1,2}(\varphi)$ is defined by (\ref{D12S}).
\end{Prop}
\pr 
Thanks to Lemma \ref{lemmaapproxG1}, we have
\begin{align*}
&   \iint_{\r^2} \hat{\rho}^N(x,y) \hat{\omega}^N(x,y) \varphi(x,y) \,dxdy  - \iint_{\r^2} \rho^0(x,y)
\omega^0(x,y) \varphi(x,y) \,dxdy  \\
= & - \iint_{\r^2\setminus C_{L}} \rho^0(x,y) \omega^0(x,y) 
\varphi(x,y) \,dxdy +  \sum_{i=0}^{N-1} \sum_{j=0}^{N-1} \int_{\overline{x}_i^N}^{\overline{x}_{i+1}^N}
\int_{\overline{y}_j^N}^{\overline{y}_{j+1}^N} \rho^0(x,y) (W_{i,j}^N -\omega^0(x,y)) \varphi(\overline{x}_i^N,\overline{y}_j^N)\,  \,dydx\\
 & +  \sum_{i=0}^{N-1} \sum_{j=0}^{N-1} W_{i,j}^N \int_{\overline{x}_i^N}^{\overline{x}_{i+1}^N}
\int_{\overline{y}_j^N}^{\overline{y}_{j+1}^N} (a_{i,j}^N-\rho^0(x,y)) ( \varphi(x,y) - \varphi(\overline{x}_i^N,\overline{y}_j^N))\,  \,dydx.
\end{align*}
Since $\rho^0$ have a compact support in $C_{L}$, we have
$\ds \iint_{\r^2\setminus C_{L}} \rho^0(x,y) \varphi(x,y) \,dxdy=0$.
From Lemma \ref{lemmaapproxsupplG1}, we have
\begin{align*}
& \left| \sum_{i=0}^{N-1} \sum_{j=0}^{N-1} \int_{\overline{x}_i^N}^{\overline{x}_{i+1}^N}
\int_{\overline{y}_j^N}^{\overline{y}_{j+1}^N} \rho^0(x,y) (W_{i,j}^N -\omega^0(x,y)) \varphi(\overline{x}_i^N,\overline{y}_j^N)\,  \,dydx \right| \\
\leq &  \frac{\max(\Delta_1,\Delta_2)}{N} (\|\partial_x \omega^0\|_\infty +
\|\partial_y \omega^0\|_\infty) \|\varphi\|_\infty \iint_{\r^2} \rho^0(x,y) \,dxdy.
\end{align*}
By an easy adaptation of Lemma \ref{lemmaapprox2}
and Lemma \ref{lemmaapprox3}, we have
\begin{align*}
& \left| \sum_{i=0}^{N-1} \sum_{j=0}^{N-1} W_{i,j}^N \int_{\overline{x}_i^N}^{\overline{x}_{i+1}^N}
\int_{\overline{y}_j^N}^{\overline{y}_{j+1}^N} (a_{i,j}^N-\rho^0(x,y)) ( \varphi(x,y) - \varphi(\overline{x}_i^N,\overline{y}_j^N))\,  \,dydx \right|  \\
\leq &  \|\omega^0\|_\infty  ( \|\partial_x \varphi\|_\infty+ \|\partial_y \varphi\|_\infty) \frac{\max(\Delta_1,\Delta_2)}{N} 
 \iint_{\r^2} \rho^0(x,y) \,dxdy \\
 & + \|\omega^0\|_\infty \|\rho^0\|_\infty  ( \|\partial_x \varphi\|_\infty+ \|\partial_y \varphi\|_\infty) \frac{\Delta_1 \Delta_2 \max(\Delta_1,\Delta_2)}{N}
\end{align*}
and we get the result.
\CQFD

This proves the remaining part of Theorem \ref{TH3}.
In the case 2 of the initial data, that is to say where $\rho^0 \in L^1(\r^2)$ doesn't not a compact support,
we have the following result.
\begin{Prop}  \label{CV_L1_SNC-G}
Let $\rho^0 \in L^1(\r^2)\cap L^{\infty}(\r^2)$ such that $\rho^0 \geq 0$ and let $\omega^0 \in L^\infty(\r^2)$.
Let $\eps >0$.
Consider $L>0$ such that (\ref{hypotheseepsL}).
We assume furthermore that $\partial_x \omega^0, \partial_y \omega^0 \in L^\infty(\r^2)$.
Let $N\in \N^*$ and
$\hat{\rho}^{N,\eps}, \hat{\rho}^{N,\eps} \hat{\omega}^{N,\eps}$ defined by (\ref{densitediscrete}) and (\ref{quantitydiscrete}) with 
(\ref{defpreli1})-(\ref{defpreli2}) and (\ref{initialdiscretquantity}).
Then for any function $\varphi \in C(\r^2)\cap L^\infty(\r^2)$ such that $\partial_x \varphi, \partial_y \varphi \in L^\infty(\r^2)$,
we have
\begin{align} \nonumber
& \left| \iint_{\r^2} \hat{\rho}^{N,\eps}(x,y) \hat{\omega}^{N,\eps}(x,y) \varphi(x,y) \,dxdy  -
 \iint_{\r^2} \rho^0(x,y) \omega^0(x,y) \varphi(x,y) \,dxdy \right| \\
& \qquad\qquad\qquad \leq    \eps  \|\varphi \|_\infty \|\omega^0\|_\infty + \frac{K_\eps(\varphi)+D_\eps(\varphi)\|\omega^0\|_\infty}{N}.
 \label{Estimation4bS}
\end{align}
where $K_\eps(\varphi)$ is defined by (\ref{KepsS}) and $D_\eps(\varphi)$ is defined by (\ref{DepsS}).
\end{Prop}
\pr  
Thanks to Lemma \ref{lemmaapproxG1}, we have
\begin{align*}
&   \iint_{\r^2} \hat{\rho}^{N,\eps}(x,y) \hat{\omega}^{N,\eps}(x,y) \varphi(x,y) \,dxdy  - \iint_{\r^2} \rho^0(x,y)
\omega^0(x,y) \varphi(x,y) \,dxdy  \\
= & - \iint_{\r^2\setminus C_{L}} \rho^0(x,y) \omega^0(x,y) 
\varphi(x,y) \,dxdy +  \sum_{i=0}^{N-1} \sum_{j=0}^{N-1} \int_{\overline{x}_i^N}^{\overline{x}_{i+1}^N}
\int_{\overline{y}_j^N}^{\overline{y}_{j+1}^N} \rho^0(x,y) (W_{i,j}^N -\omega^0(x,y)) \varphi(\overline{x}_i^N,\overline{y}_j^N)\,  \,dydx\\
 & +  \sum_{i=0}^{N-1} \sum_{j=0}^{N-1} W_{i,j}^N \int_{\overline{x}_i^N}^{\overline{x}_{i+1}^N}
\int_{\overline{y}_j^N}^{\overline{y}_{j+1}^N} (a_{i,j}^N-\rho^0(x,y)) ( \varphi(x,y) - \varphi(\overline{x}_i^N,\overline{y}_j^N))\,  \,dydx.
\end{align*}
First we have
\begin{align*}
\left| \iint_{\r^2\setminus C_L} \rho^0(x,y) \varphi(x,y) \omega^0(x,y) \,dxdy \right| & \leq \|
 \varphi \|_\infty \|\omega^0 \|_\infty  \iint_{\r^2\setminus C_L} \rho^0(x,y)  \,dxdy
\leq  \| \varphi \|_\infty \|\omega^0 \|_\infty \eps.
\end{align*}
From Lemma \ref{lemmaapproxsupplG1}, we have
\begin{align*}
& \left| \sum_{i=0}^{N-1} \sum_{j=0}^{N-1} \int_{\overline{x}_i^N}^{\overline{x}_{i+1}^N}
\int_{\overline{y}_j^N}^{\overline{y}_{j+1}^N} \rho^0(x,y) (W_{i,j}^N -\omega^0(x,y)) \varphi(\overline{x}_i^N,\overline{y}_j^N)\,  \,dydx \right| \\
\leq &  \frac{\Delta}{N} (\|\partial_x \omega^0\|_\infty +
\|\partial_y \omega^0\|_\infty) \|\varphi\|_\infty \iint_{\r^2} \rho^0(x,y) \,dxdy.
\end{align*}
By an adaptation of Lemma \ref{lemmaapprox2}
and Lemma \ref{lemmaapprox3} and since $\Delta_2=\Delta_1=\Delta=2L$, we have
\begin{align*}
& \left| \sum_{i=0}^{N-1} \sum_{j=0}^{N-1} W_{i,j}^N \int_{\overline{x}_i^N}^{\overline{x}_{i+1}^N}
\int_{\overline{y}_j^N}^{\overline{y}_{j+1}^N} (a_{i,j}^N-\rho^0(x,y)) ( \varphi(x,y) - \varphi(\overline{x}_i^N,\overline{y}_j^N))\,  \,dydx \right|  \\
\leq &  \|\omega^0\|_\infty  ( \|\partial_x \varphi\|_\infty+ \|\partial_y \varphi\|_\infty) \frac{2L}{N} 
 \iint_{\r^2} \rho^0(x,y) \,dxdy \\
 & + \|\omega^0\|_\infty \|\rho^0\|_\infty  ( \|\partial_x \varphi\|_\infty+ \|\partial_y \varphi\|_\infty) \frac{8L^3}{N}
\end{align*}
and we get the result.
\CQFD

This proves the remaining part of Theorem \ref{TH4}.


\end{document}